\documentclass[11pt]{amsart}
\usepackage{amsmath,amstext,amssymb,amsopn,amsthm,mathrsfs}
\usepackage[top=4cm, bottom=3cm, left=3cm, right=3cm]{geometry}
\usepackage[shortlabels]{enumitem}
\usepackage{amsfonts}
\usepackage{amsthm}
\usepackage[T1]{fontenc}
\usepackage[makeroom]{cancel}
\usepackage{color}
\usepackage[dvips]{graphicx}

\newtheorem{thm}{Theorem}[section]
\newtheorem{df}{Definition}[section]
\newtheorem{cor}{Corollary}[section]
\newtheorem{lem}{Lemma}[section]
\newtheorem{prob}{Proposition}[section]

\numberwithin{equation}{section}

\begin{document}
\title[Notes on Harmonic Analysis]{Notes on Harmonic Analysis\\
 Part I: The
Fourier Transform}
\author{Kecheng Zhou}
\address{Kecheng Zhou, Department of Mathematics, California State
University, Sacramento, CA, 95819, USA.}
\author{M. Vali Siadat}
\address{M. Vali Siadat, Department of Mathematics and Statistics, Loyola
University Chicago, Chicago, IL, 60660, USA}

\begin{abstract}
Fourier Transforms is a first in a series of monographs we present on
harmonic analysis. Harmonic analysis is one of the most fascinating areas of
research in mathematics. Its centrality in the development of many areas of
mathematics such as partial differential equations and integration theory
and its many and diverse applications in sciences and engineering fields
makes it an attractive field of study and research.

The purpose of these notes is to introduce the basic ideas and theorems of
the subject to students of mathematics, physics or engineering sciences. Our
goal is to illustrate the topics with utmost clarity and accuracy, readily
understandable by the students or interested readers. Rather than providing
just the outlines or sketches of the proofs, we have actually provided the
complete proofs of all theorems. This will illuminate the necessary steps
taken and the machinery used to complete each proof.

The prerequisite for understanding the topics presented is the knowledge of
Lebesgue measure and integral. This will provide ample mathematical
background for an advanced undergraduate or a graduate student in
mathematics.
\end{abstract}

\maketitle

\section{Fourier Transforms for $L^1( \mathbb{R})$}

\begin{df}
For $f\in L^{1}(\mathbb{R}),$ the Fourier transform $\hat{f}$ of $f$ is
defined as 
\begin{equation}
\hat{f}(y)=\int_{-\infty }^{\infty }f(x)e^{-ixy}dx
\end{equation}

for all real $y\in \mathbb{R}.$
\end{df}

It is easy to see that Fourier transform is a lineaer operator, i.e., $(f+g)%
\hat{(}y)=\hat{f}(y)+\hat{g}(y)$ and $(kf)\hat{(}y)=k\hat{f}(y).$ Also using
simple integration techniques it can easily be shown that $\ \hat{f}%
(y+t)=e^{ity}\hat{f}(y)$ and $\ \hat{f}(ky)=\frac{1}{k}\hat{f}(\frac{y}{k}).$

\begin{thm}
If $f\in L^{1}(\mathbb{R}),$ then $\hat{f}(y)$ is uniformly continuous and
bounded in $\mathbb{R}.$
\end{thm}

\textbf{Proof:}\quad Clearly, $|\hat{f}(y)|\leq ||f||_{1}$ for all $y.$
Moreover, 
\begin{eqnarray*}
|\hat{f}(y+h)-\hat{f}(y)|&=&|\int_{-\infty }^{\infty}
f(x)(e^{-ix(y+h)}-e^{-ixy})dx| \\
&\leq& \int_{-\infty }^{\infty }|f(x)||e^{-ixh}-1|dx.
\end{eqnarray*}%
The integrand on the right side converges to $0$ as $h\rightarrow 0$ and is
dominated by $2|f(x)|\in L^{1}(\mathbb{R}).$ So, by Lebesgue's dominated
convergence theorem, $\hat{f}$ is uniformly continuous.\qed

\begin{thm}[Riemann-Lebesgue Lemma]
If $f\in L^{1}(\mathbb{R}),$ then $\hat{f}(y)\rightarrow 0$ as $y\rightarrow
\pm \infty . $
\end{thm}

\textbf{Proof:}\quad First suppose that $f$ is a characteristic function of
an interval $[a,b]$. Its Fourier transform is 
\begin{equation*}
\int_{a}^{b}e^{-ixy}dx=\frac{e^{-iay}-e^{-iby}}{iy},\quad y\neq 0,
\end{equation*}%
which tends to zero. Therefore, a linear combination of characteristic
functions of intervals, i.e., a step function, satisfies the
Riemann-Lebesgue lemma. Such functions are also dense in $L^{1}(\mathbb{R}).$
Now let $f\in L^{1}(\mathbb{R})$ and let $f_{n}\in L^{1}(\mathbb{R})$ be a
sequence of step functions such that $f_{n}\rightarrow f$ in $L^{1}(\mathbb{R%
}).$ Then 
\begin{equation*}
|\hat{f_{n}}(y)-\hat{f}(y)|=|(f_{n}-f)\hat{\,}(y)|\leq
||f_{n}-f||_{1}\rightarrow 0.
\end{equation*}%
Note that the limit is uniform in $y\in \mathbb{R}.$ Since 
\begin{equation*}
|\hat{f}(y)|\leq |\hat{f_{n}}(y)-\hat{f}(y)|+|\hat{f_{n}}(y)|,
\end{equation*}%
we can choose $n$ large enough so that the first term on the right is small
and then for that fixed $n,$ we let $|y|$ large enough so that the second
term is also small. This completes the proof. \qed

\begin{thm}
Suppose that $f(x)(1+|x|)$ is integrable. Then, 
\begin{equation}
(\hat{f})^{\prime }(y)=(-ixf(x))\hat{\,}(y).
\end{equation}
\end{thm}

\textbf{Proof:}\quad Note that, by assumption, both $f$ and $xf(x)$ are
integrable. We write 
\begin{eqnarray*}
(\hat{f})^{\prime }(y) &=&\lim_{h\rightarrow 0}\int_{-\infty }^{\infty }f(x)%
\frac{e^{-ix(y+h)}-e^{-ixy}}{h}dx \\
&=&\lim_{h\rightarrow 0}\int_{-\infty }^{\infty }f(x)e^{-ixy}\frac{e^{-ixh}-1%
}{h}dx.
\end{eqnarray*}%
Note that the integrand converges to $f(x)e^{-ixy}(-ix)$ pointwise as $%
h\rightarrow 0$ and $|f(x)e^{-ixy}\frac{e^{-ixh}-1}{h}|\leq |xf(x)|$ for all
small $|h|.$ \footnote{%
Estimating the remainder (both Lagrange form and integral form) of Taylor's
series for $e^{ix}$ we obtain the estimation 
\begin{equation*}
|e^{ix}-\sum_{k=0}^{n}\frac{(ix)^{k}}{k!}|\leq \min (\frac{|x|^{n+1}}{(n+1)!}%
,\frac{2|x|^{n}}{n!}).
\end{equation*}%
Note that the first estimate is better for small $|x|,$ while the second is
better for large $|x|.$ Choosing $n=0$ and considering small $|h|$ we get
the inequality in context.} Hence, by Lebesgue's dominated convergence
theorem , 
\begin{equation*}
(\hat{f})^{\prime }(y)=\lim_{h\rightarrow 0}\int_{-\infty }^{\infty
}f(x)e^{-ixy}\frac{e^{-ixh}-1}{h}dx=\int_{-\infty }^{\infty
}(-ixf(x))e^{-ixy}dx. \qed
\end{equation*}

\begin{thm}
If $f$ is continuously differentiable with compact support, then 
\begin{equation}
(f^{\prime})\hat(y)=iy\hat{f}(y).
\end{equation}
\end{thm}

\textbf{Proof:}\quad Integration by parts.\qed

\begin{df}
The convolution of $f$ and $g$ is defined as 
\begin{equation}
(f\ast g)\,(x)=\int_{-\infty }^{\infty }f(x-t)g(t)dt,
\end{equation}
whenever the integral exists.\newline
\end{df}

In the following, $C(\mathbb{R})$ denotes the space of all continuous
functions on $\mathbb{R}$ with $||f||_C=sup_{x\in \mathbb{R}}|f(x)|<\infty$
and $C_0(\mathbb{R})$ the space of all continuous functions on $\mathbb{R}$
that vanishes at infinity, i.e., for any $\epsilon>0,$ there is a compact $%
F\subset \mathbb{R}$ such that $|f(x)|<\epsilon$ for $x\not\in F.$ Then by
F. Riesz' theorem, $(C_0(\mathbb{R}))^*=M(\mathbb{R}),$ where $M(\mathbb{R})$
is the space of complex regular Borel measures on $R.$\footnote{%
Every complex measure is bounded, see Hewitt and Stromberg \cite%
{Hewitt&Stromberg} .} Since $C_0(\mathbb{R})$ is separable (continuous
functions with compact support are dense in $C_0(\mathbb{R})$), every
bounded subset of $M(\mathbb{R})$ is weak* sequentially compact. Note that $%
L^1(\mathbb{R})$ is contained in $M(\mathbb{R}),$ if we identify $f\in L^1(%
\mathbb{R})$ with the measure $f(x)dx.$\newline

\begin{thm}
Let $f\in L^{p}(\mathbb{R}),$ $1\leq p\leq \infty ,$ and $g\in L^{q},\frac{1%
}{p}+\frac{1}{q}=1.$ Then $(f\ast g)(x)$ exists everywhere, belongs to $C(%
\mathbb{R}),$ and $||f\ast g||_{C}\leq ||f||_{p}||g||_{q}.$ Moreover, if $%
1<p<\infty ,$ or if $p=1$ and $g\in C_{0}(\mathbb{R}),$ then $f\ast g\in
C_{0}(\mathbb{R}),$ i.e., $f\ast g\in C(\mathbb{R})$ and $%
\lim_{|x|\rightarrow \infty }|(f\ast g)(x)|=0.$
\end{thm}

\textbf{Proof:}\quad Let $1\leq p<\infty .$ By H\"{o}lder's inequality, $%
|(f\ast g)(x)|\leq ||f||_{p}||g||_{q}$ and so $(f\ast g)(x)$ exists for
every $x\in \mathbb{R}.$ Furthermore, 
\begin{equation*}
|(f\ast g)(x+h)-(f\ast g)(x)|\leq ||f(\cdot +h)-f(\cdot )||_{p}||g||_{q},
\end{equation*}%
and therefore $f\ast g\in C(\mathbb{R})$ by the continuity of $f$ in mean.
If $p=\infty , $ the roles of $f$ and $g$ can be interchanged.

Now let $1<p<\infty $ (obviously $1<q<\infty $ also). Given $\epsilon >0,$
there is a finite interval $[-a,a]$ such that 
\begin{equation*}
\int_{|t|\geq a}|f(t)|^{p}dt\leq \epsilon ^{p}\mbox{ and }\int_{|t|\geq
a}|g(t)|^{q}dt\leq \epsilon ^{q}.
\end{equation*}%
If $x\in \mathbb{R}$ is such that $|x|>2a,$ then $[x-a,x+a]$ is contained in 
$\{t\in \mathbb{R}:|t|>a\},$ and hence 
\begin{eqnarray*}
|(f\ast g)(x)| &\leq &(\int_{-a}^{a}+\int_{|t|\geq a})|f(x-t)g(t)|dt \\
&\leq &(\int_{-a}^{a}|f(x-t)|^{p}dt)^{1/p}||g||_{q}+||f||_{p}(\int_{|t|\geq
a}|g(t)|^{q}dt)^{1/q} \\
&\leq &(\int_{x-a}^{x+a}|f(t)|^{p}dt)^{1/p}||g||_{q}+||f||_{p}\epsilon \\
&\leq &(\int_{|t|>a}|f(t)|^{p}dt)^{1/p}||g||_{q}+||f||_{p}\epsilon \leq
\epsilon (||g||_{q}+||f||_{p}).
\end{eqnarray*}%
Thus, $(f\ast g)(x)$ tends to $0$ as $|x|\rightarrow\infty,$ giving $f\ast
g\in C_{0}(\mathbb{R}).$ The same method of proof applies for the case $p=1,$
$g\in C_{0}(\mathbb{R}).$ \qed

\begin{thm}
If $f,g\in L^{1}(\mathbb{R}),$ then $(f\ast g)(x)$ exists a.e. and $||f\ast
g||_{1}\leq ||f||_{1}||g||_{1}.$ Moreover, 
\begin{equation}
(f\ast g)\hat{\,}(y)=\hat{f}(y)\cdot \hat{g}(y).
\end{equation}
\end{thm}

\textbf{Proof:}\quad Note that the integral ${\displaystyle\int
|f(x-t)g(t)|dx}$ exists for a.e. $t$ and 
\begin{equation*}
\int |f(x-t)g(t)|dx=|g(t)|\cdot ||f||_{1}.
\end{equation*}%
Also note that the expression on the right belongs to $L^{1}(\mathbb{R}).$
Hence, the integral 
\begin{equation*}
\int (\int |f(x-t)g(t)|dx)dt=||f||_{1}||g||_{1}
\end{equation*}%
exists as a finite number. Therefore, by Fubini's theorem the integral 
\begin{equation*}
\int (\int |f(x-t)g(t)|dt)dx
\end{equation*}%
exists and is equal to $||f||_{1}||g||_{1}.$ This implies that ${%
\displaystyle\int |f(x-t)g(t)|dt}$ exists a.e. and belongs to $L^{1}.$

To prove $(f\ast g)\hat{\,}(y)=f\hat{\,}(y)\cdot g\hat{\,}(y),$ we observe
that 
\begin{eqnarray*}
(f\ast g)\hat{\,}(y) &=&\int (\int f(x-t)g(t)dt)e^{-ixy}dx \\
&=&\int g(t)e^{-ity}(\int f(x-t)e^{-iy(x-t)}dx)dt=\hat{f}(y)\hat{g}(y).
\end{eqnarray*}%
The change in the order of integration is justified by Fubini's theorem. %
\qed \newline

It is easy to see that convolution obeys the commutative and distributive
laws of algebra in $L^{1}(\mathbb{R}),$ i.e., $f\ast g=g\ast f$ and $f\ast
(g+h)=f\ast g+f\ast h.$ The natural question is whether there is a
multiplicative identity, i.e., given $f\in L^{1}(\mathbb{R}),$ is there $%
e\in L^1(\mathbb{R})$ such that $f\ast e=f?$ The answer is, in general, no
since convolution exhibits continuity property and cannot be equal to a
discontinuous $f.$ However, we may seek a sequence of functions $e_{n},$
called approximate identity, with the property that $e_{n}\ast f\rightarrow
f.$\newline

\begin{df}
An approximate identity $e_{n}$ on $\mathbb{R}$ is a sequence of functions $%
e_{n}$ such that $e_{n}\geq 0,$ \,\,\,${\displaystyle\int e_{n}(x)dx=1},$
and \,\,\,for each $\delta >0,$ 
\begin{equation*}
\lim_{n\rightarrow \infty} \int_{|x|>\delta}e_{n}(x)dx=0.
\end{equation*}
\end{df}

\begin{thm}
If $f\in C_{0}(\mathbb{R}),$ then $e_{n}\ast f\rightarrow f$ uniformly. If $%
f\in L^{p}(\mathbb{R}),$ $1\leq p<\infty ,$ then $e_{n}\ast f\rightarrow f$
in $L^{p}(\mathbb{R}).$ If $f\in L^{\infty }(\mathbb{R}),$ then $e_{n}\ast
f\rightarrow f$ in the weak* topology of $L^{\infty }(\mathbb{R})$ as a dual
of $L^{1}(\mathbb{R}),$ that is, ${\displaystyle\int (e_{n}\ast
f)(x)g(x)dx\rightarrow \int f(x)g(x)dx}$ for all $g\in L^{1}(\mathbb{R}).$
\end{thm}

\textbf{Proof:}\quad Note that if $f\in C_{0}(\mathbb{R}),$ then $f$ is
uniformly continuous on $\mathbb{R}$ and for any given $\epsilon >0$, there
is a $\delta >0$ such that for any $t$ with $|t|<\delta ,$ $%
|f(x-t)-f(x)|<\epsilon $ for all $x\in \mathbb{R}.$ Hence, 
\begin{eqnarray*}
&&|(e_{n}\ast f)(x)-f(x)|\leq \int_{R}|f(x-t)-f(x)|e_{n}(t)dt \\
&=&\int_{|t|<\delta }|f(x-t)-f(x)|e_{n}(t)dt+\int_{|t|\geq \delta
}|f(x-t)-f(x)|e_{n}(t)dt \\
&\leq &\epsilon +2M\int_{|t|>\delta }e_{n}(t)dt,
\end{eqnarray*}%
where $M=\sup_{x\in \mathbb{R}}|f(x)|.$ Since $\lim_{n\rightarrow \infty
}\int_{|x|>\delta }e_{n}(x)dx=0,$\thinspace\ $e_{n}\ast f\rightarrow f$
uniformly. In the case of $f\in L^{\infty }(\mathbb{R}),$ the proof is
similar.

If $f\in L^p(\mathbb{R}), 1\leq p<\infty,$ then 
\begin{eqnarray*}
\int_\mathbb{R} |(e_n*f)(x)-f(x)|^pdx&\leq& \int_\mathbb{R} |\int_\mathbb{R}
(f(x-t)-f(x))e_n(t)dt|^p dx \\
&\leq &\int_\mathbb{R} (\int_\mathbb{R} |f(x-t)-f(x)|^p\,|e_n(t)|dx) dt \\
&=&\int_\mathbb{R} ||f(\cdot -t)-f(\cdot )||_p^p \, |e_n(t)| dt.
\end{eqnarray*}
Given any $\epsilon>0,$ there is a $\delta>0$ such that $||f(\cdot
-t)-f(\cdot )||_p<\epsilon$ whenever $|t|<\delta.$ Hence, 
\begin{eqnarray*}
&& \int_R ||f(\cdot -t)-f(\cdot )||_p^p \, |e_n(t)| dt \\
&=&\int_{|t|<\delta}||f(\cdot -t)-f(\cdot )||_p^p \, e_n(t)dt+\int_{|t|\geq
\delta} ||f(\cdot -t)-f(\cdot )||_p^p \, e_n(t)dt \\
&\leq & \epsilon^p+2||f||_p^p\int_{|t|>\delta}e_n(t)dt.
\end{eqnarray*}
Since $\lim_{n\rightarrow \infty}\int_{|x|>\delta}e_{n}(x)dx=0,$\, the
result follows. \qed

\begin{thm}
If $f$ has compact support and a continuous derivative, and $g\in L^{1}(%
\mathbb{R}),$ then $f\ast g\in L^{1}(\mathbb{R})$ has a continuous
derivative.
\end{thm}

\textbf{Proof:}\quad First, we prove 
\begin{equation*}
\frac{d}{dx}(\int f(x-t)g(t)dt)=\int f^{\prime }(x-t)g(t)dt,
\end{equation*}
which is showing that 
\begin{equation*}
\lim_{h\rightarrow 0}\int_{-\infty }^{\infty }(\frac{f(x+h-t)-f(x-t)}{h}
)g(t)dt=\int f^{\prime }(x-t)g(t)dt.
\end{equation*}
Note that the integrand on the left converges to $f^{\prime }(x-t)g(t)$
pointwise (in $t$) as $h\rightarrow 0.$ Moreover, $\frac{f(x+h-t)-f(x-t)}{h}
=f^{\prime }(c),$ where $c$ is between $x+h-t$ and $x-t.$ If $f$ has compact
support $S,$ then so does $f^{\prime }.$ Therefore, $|\frac{f(x+h-t)-f(x-t)}{
h}|=|f^{\prime }(c)|\leq sup_{c\in S}|f^{\prime }(c)|\leq M$ with some $M>0$
for all $t\in (-\infty ,\infty ).$ Now the desired limit follows from
Lebesgue's dominated convergence theorem.

To prove that $\int f^{\prime}(x-t)g(t) dt$ is continuous, we note that 
\begin{eqnarray*}
&&|\int f^{\prime}(x+h-t)g(t)dt-\int f^{\prime}(x-t)g(t)| \\
&=& |\int f^{\prime}(t)(g(x+h-t)-g(x-t))dt| \\
&\leq& ||f^{\prime}||_c ||g(\cdot+h)-g(\cdot)||_1.
\end{eqnarray*}
Then the (uniform) continuity of ${\displaystyle\int f^{\prime}(x-t)g(t)dt}$
follows from the continuity of $g$ in mean. \qed\newline

The following corollary follows immediately from Theorems 1.7 and 1.8.

\begin{cor}
Let $e_n(x)$ be an approximate identity with compact support and continuous
derivative. Then for any $f\in L^1(\mathbb{R}),$ $e_n*f$ provides a
continuously differentiable approximation to $f$ in $L^1(\mathbb{R}).$
\end{cor}

\textbf{Proof:} An obvious result from Theorem 1.8. \qed

\begin{thm}
Let $\phi(x)\geq 0$ be a function defined on $\mathbb{R}$ such that $\phi$
has compact support and continuous derivative, and $\int \phi(x)dx=1.$ Then $%
e_n(x)=n\phi(nx)$ is an approximate identity with compact support and
continuous derivative.
\end{thm}

\textbf{Proof:}\quad We only need to show that for each $\epsilon >0,$ ${%
\displaystyle\int_{|t|\geq \epsilon }e_{n}(t)dt=0.}$ In fact, 
\begin{equation*}
\int_{|t|\geq \epsilon }n\phi (nt)dt=\int_{|u|\geq n\epsilon }\phi
(u)du\rightarrow 0,\,\,
\end{equation*}
as $n\rightarrow \infty.$ \qed

\begin{thm}
If $f$ and ${\displaystyle\frac{f(x)}{x}}$ \footnote{%
The assumption that ${\displaystyle\frac{f(x)}{x}}\in L^{1}(\mathbb{R})$
simply emphasizes that $f$ behaves like a \emph{positive} power of $x$ at $%
x=0.$ For example, If $f\in Lip(\alpha )$ for $0<\alpha \leq 1$ at $x=0,$
and $f(0)=0,$ then ${\displaystyle\frac{f(x)}{x}}\in L^{1}(\mathbb{R}).$}
are both integrable, then 
\begin{equation*}
\lim_{A,B\rightarrow \infty }\int_{-B}^{A}\hat{f}(y)dy=0.
\end{equation*}
\end{thm}

\textbf{Proof:}\quad Observe that 
\begin{equation*}
\int_{-B}^{A}\hat{f}(y)dy=\int_{-B}^{A}\int_{-\infty }^{\infty
}f(x)e^{-ixy}dxdy.
\end{equation*}%
The integrand on the right side is integrable over the product space, so, by
Fubini's theorem, 
\begin{eqnarray*}
\int_{-B}^{A}\int_{-\infty }^{\infty }f(x)e^{-ixy}dxdy&=&\int_{-\infty
}^{\infty }\int_{-B}^{A}f(x)e^{-ixy}dydx \\
&=&\int_{-\infty }^{\infty }f(x)\frac{e^{iBx}-e^{-iAx}}{ix}dx.
\end{eqnarray*}
The last integral tends $0$ $($as $A,B\rightarrow \infty ),$ by
Riemann-Lebesgue lemma. \qed\newline

To derive the following inversion theorem, we need a simple fact, which can
be verified by a straightforward calculation: If $g(x)=e^{-|x|},$ then $\hat{%
g}(y)={\ \displaystyle\frac{2}{1+y^{2}}}$ and ${\displaystyle\frac{1}{2\pi }
\int_{-\infty }^{\infty }\hat{g}(y)dy=1.}$

\begin{cor}
If $f$ is integrable in $\mathbb{R}$ and satisfies a Lipschitz condition at $%
t,$ then 
\begin{equation*}
f(t)=\lim_{A,B\rightarrow \infty }\frac{1}{2\pi }\int_{-B}^{A}\hat{f}%
(y)e^{ity}dy.
\end{equation*}
\end{cor}

That is, $f$ is the inverse Fourier transform of $\hat{f}{.}$\newline

\textbf{Proof:}\quad If $t\neq 0,$ let $h(x)=f(x+t).$ If $f(x)$ satisfies $%
|f(x)-f(t)|\leq K|x-t|^{\alpha }$ for $x$ near $t$ then for $x$ near $0$, $%
h(x)$ satisfies the Lipschitz condition at $t=0$: $|h(x)-h(0)|=|f(x+t)-f(t)|%
\leq K|x|^{\alpha }.$ Therefore, if we can show the corollary for $t=0,$
then for $t\neq 0,$ 
\begin{equation*}
f(t)=h(0)=\lim_{A,B\rightarrow \infty }\frac{1}{2\pi }\int_{-B}^{A}\hat{h}%
(y)dy=\lim_{A,B\rightarrow \infty }\frac{1}{2\pi }\int_{-B}^{A}\hat{f}%
(y)e^{ity}dy.\,\,
\end{equation*}

We may now assume that $t=0.$ Since $f$ satisfies the Lipschitz condition at 
$0,$ if $f(0)=0,$ it follows that ${\displaystyle\frac{f(x)}{x}}\in L^{1}(%
\mathbb{R}).$ Then by Theorem 1.10, ${\displaystyle\lim_{A,B\rightarrow
\infty }\int_{-B}^{A}\hat{f}(y)dy=0},$ which shows that 
\begin{equation*}
f(t)=\lim_{A,B\rightarrow \infty }\frac{1}{2\pi }\int_{-B}^{A}\hat{f}%
(y)e^{ity}dy
\end{equation*}%
holds as $t=0$ and $f(0)=0.$

If $f(0)\neq 0,$ we reduce it to the case $f(0)=0.$ Let $\phi
(x)=f(x)-f(0)g(x),$ where $g(x)=e^{-|x|}.$ Then $\phi (0)=0$ and $\phi (x)$
satisfies the Lipschitz condition at $t=0.$ Therefore, 
\begin{equation*}
\frac{1}{2\pi }\int_{-B}^{A}\hat{\phi}(y)dy\rightarrow 0,
\end{equation*}%
as $A,B\rightarrow \infty .$ That is, 
\begin{equation*}
\frac{1}{2\pi }\int_{-B}^{A}\hat{f}(y)dy-\frac{1}{2\pi }\int_{-B}^{A}f(0)%
\hat{g}(y)dy\rightarrow 0.
\end{equation*}%
It follows that 
\begin{equation*}
\frac{1}{2\pi }\int_{-B}^{A}\hat{f}(y)dy\rightarrow f(0). \qed
\end{equation*}%
\textbf{Remark:} For $f\in L^{1}(\mathbb{R}),$ $\hat{f}$ need not be in $%
L^{1}(\mathbb{R}).$ Therefore, the above integral ${\displaystyle%
\int_{-B}^{A}\hat{f}(y)e^{ity}dy}$ has to be understood as the limit of the
integral from $-B$ to $A$ as $A,B\rightarrow \infty .$ Note that a function $%
f\geq 0$ is integrable on $\mathbb{R}$ if $\lim_{A,B\rightarrow \infty
}\int_{-B}^{A}f(x)dx$ exists. Therefore, a function $f$ is integrable on $%
\mathbb{R}$ if both $f^{+}$ and $f^{-}$ are integrable on $\mathbb{R}$. By
this definition, $f$ and $|f|$ are either both integrable or not integrable.
Hence, it may happen that $\hat{f}\not\in L^{1}(\mathbb{R}),$ yet the above
limit exists. Let $m(x)=\chi _{\lbrack -a,a]}(x).$ Clearly, $m\hat{\,}(y)={%
\displaystyle\frac{2\sin ay}{y}}\not\in L^{1}(\mathbb{R}),$ but 
\begin{equation*}
{\displaystyle m(x)=\lim_{A,B\rightarrow \infty }\frac{1}{2\pi }\int_{-B}^{A}%
\hat{m}(y)e^{ixy}dy,\,\,\,\,\,x\neq \pm 1.}
\end{equation*}
With improper Riemann integral in mind, we may say $f$ equals the inverse
Fourier transform of $\hat{f}$ at each Lipschitz point.\newline

\begin{thm}
\begin{equation}
\frac{1}{2\pi}\int_{-\infty}^{\infty}\frac{2}{1+t^2}e^{ixt}dt =e^{-|x|}.
\end{equation}
\end{thm}

\textbf{Proof:}\quad Assuming that $x>0$ and integrating 
\begin{equation*}
I=\frac{1}{2\pi }\int_{\Gamma _{R}}\frac{2}{1+z^{2}}e^{ixz}dz,
\end{equation*}%
where $\Gamma _{R}$ consists of the upper semicircle $C_{R}$ and the line
segment $[-R,R]$ on the $x$-axis, we see that 
\begin{equation*}
I=Res_{z=i}\frac{2}{1+z^{2}}e^{ixz}=e^{-x}
\end{equation*}%
and that the integral along $[-R,R]$ gives 
\begin{equation*}
I=\frac{1}{2\pi }\int_{\Gamma _{R}}\frac{2}{1+t^{2}}e^{ixt}dt,
\end{equation*}%
while, if $x>0,$ then 
\begin{equation*}
\frac{1}{2\pi }\int_{C_{R}}\frac{2}{1+z^{2}}e^{ixz}dz\rightarrow 0
\end{equation*}%
as $R\rightarrow \infty .$ \footnote{%
(Jordan's Lemma) Suppose that $f$ is an analytic function in the upper half
plane except at a finite number of singularities and $|f(z)|\rightarrow 0$
as $|z|\rightarrow \infty $ for $0\leq Arg(z)\leq \pi $. Then, if $x>0,$ $%
\int_{C_{R}}e^{ixz}f(z)dz\rightarrow 0$ as $R\rightarrow \infty .$} Thus, if 
$x>0,$ $\frac{1}{2\pi }\int_{-\infty }^{\infty }\frac{2}{1+t^{2}}%
e^{ixt}dt=e^{-x}.$ Similarly, if $x<0$ then $\frac{1}{2\pi }\int_{-\infty
}^{\infty }\frac{2}{1+t^{2}}e^{ixt}dt=e^{x}.$ Hence, 
\begin{equation*}
\frac{1}{2\pi }\int_{-\infty }^{\infty }\frac{2}{1+t^{2}}e^{ixt}dt=e^{-|x|}. %
\qed
\end{equation*}

In the following, we will calculate the Fourier transform of a Gaussian
function which will be useful in proving the inversion theorem. The theorem
below simply states that Fourier transform of a Gaussian function is a
Gaussian.

\begin{thm}
\begin{equation}
(e^{-x^{2}})\hat{\,}(y)=\sqrt{\pi }e^{-y^{2}/4}.
\end{equation}
\end{thm}

\textbf{Proof:}\quad First, let $u$ be real. We have 
\begin{equation*}
\int_{-\infty}^\infty e^{-x^2+2xu}dx=e^{u^2} \int_{-\infty}^\infty
e^{-(x-u)^2}dx =e^{u^2}\int_{-\infty}^\infty e^{-t^2}dt=\sqrt{\pi}e^{u^2}.
\end{equation*}

Clearly, the function defined by 
\begin{equation*}
\int_{-\infty }^{\infty }e^{-x^{2}+2xz}dx
\end{equation*}%
is an entire function\footnote{%
See the theorem in complex analysis. Suppose that $f(z,w)$ is a continuous
function of $z\in D$ and $w\in C,$ where $D$ is a region and $C$ is a
contour that is a piecewise smooth curve $w(t)=u(t)+iv(t),$ $t_{0}\leq t\leq
t_{1},$ with continuous $u^{\prime }$ and $v^{\prime }.$ Suppose that for
each $w\in C,$ $f(z,w)$ is an analytic function in $z\in D.$ Then $%
F(z)=\int_{C}f(z,w)dw$ is analytic in $D$ and $F^{\prime }(z)$ can be found
by differentiating under the integral sign.
\par
If $C$ is a contour going to infinity such that any bounded part of it is
regular (no sharp corner) and if the above conditions are satisfied on any
bounded part of $C,$ and if $\int_{C}f(z,w)dw$ converges uniformly in $z\in
D $, then the above results hold.}, and by above calculation, it coincides
with the entire function $\sqrt{\pi }e^{z^{2}}$ along the $x$-axis.
Therefore, for all $z$, 
\begin{equation*}
\int_{-\infty }^{\infty }e^{-x^{2}+2xz}dx=\sqrt{\pi }e^{z^{2}}.
\end{equation*}%
In particular, let $z=\frac{-iy}{2}.$ Then we have 
\begin{equation*}
\int_{-\infty }^{\infty }e^{-x^{2}-ixy}dx=\sqrt{\pi }e^{-y^{2}/4}. \qed 
\newline
\end{equation*}

\begin{thm}[Inversion Theorem]
\ Let $f\in L^1(\mathbb{R}),$ and $\hat{f}\in L^1(\mathbb{R}),$ then%
\begin{equation}
f(x)=\frac{1}{2\pi }\int_{-\infty }^{\infty }\hat{f}(y)e^{ixy}dy
\end{equation}
for almost all real $x\in \mathbb{R}$. The integral is commonly known as the
inverse Fourier transform.
\end{thm}

\textbf{Proof:}\quad Consider the Gauss-Weierstrass Kernel, $W(x,\alpha )=%
\frac{1}{\sqrt{\pi \alpha }}e^{-\frac{x^{2}}{\alpha }}.$ A straightforward
calculation shows that $W(\cdot ,\alpha )\hat{\,}(t)=e^{-\frac{\alpha t^{2}}{%
4}}.$ By integrating $\hat{f}$ against $W\hat{\,}$, and then applying
Fubini's theorem and the fact that $W(x,\alpha )$ is an approximate
identity, we get

\begin{eqnarray*}
&&\int_{-\infty }^{\infty }\hat{f}(\xi )e^{i\xi x}W(\cdot ,\alpha )\hat{\,}%
(\xi )d\xi \\
&=&\int_{-\infty }^{\infty }(\int_{-\infty }^{\infty }f(t)e^{-i\xi
t}dt)e^{i\xi x}e^{-\frac{\alpha \xi ^{2}}{4}}d\xi \\
&=&\int_{-\infty }^{\infty }f(t)(\int_{-\infty }^{\infty }e^{-\frac{\alpha
\xi ^{2}}{4}}e^{-i\xi (t-x)}d\xi )dt \\
&=&\int_{-\infty }^{\infty }f(t)2\pi W(t-x,\alpha )dt \\
&=&2\pi \int_{-\infty }^{\infty }f(x-t)W(t,\alpha )dt\rightarrow 2\pi
f(x)\,\,\,a.e.\mbox{
\ as \ }\alpha \rightarrow 0^{+}
\end{eqnarray*}%
On the other hand, by Lebesgue's dominated convergence theorem, 
\begin{equation*}
\lim_{\alpha \rightarrow 0^{+}}\int_{-\infty }^{\infty }\hat{f}(\xi )e^{i\xi
x}W(\cdot ,\alpha )\hat{\,}(\xi )d\xi =\int_{-\infty }^{\infty }\hat{f}(\xi
)e^{i\xi x}d\xi .
\end{equation*}%
The theorem follows. \qed \newline

As an application of the inversion theorem, we now prove that the Fourier
transform of a product is the convolution of the Fourier transforms.

\begin{thm}
Assume that $f,g\in L^{1}(\mathbb{R})$ and $\hat{f}\in L^{1}(\mathbb{R})$ $($%
or $\hat{g}\in L^{1}(\mathbb{R})).$ Then, 
\begin{equation}
(fg)\hat{\,}(x)=\frac{1}{2\pi }(\hat{f}\ast \hat{g})(x).
\end{equation}
\end{thm}

\textbf{Proof:}\quad By the inversion theorem, $f$ is bounded and so, $fg\in
L^{1}(\mathbb{R}).$ Hence,%
\begin{eqnarray*}
(fg)\hat{\,}(x) &=&\int_{-\infty }^{\infty }f(y)g(y)e^{-ixy}dy \\
&=&\int_{-\infty }^{\infty }g(y)e^{-ixy}(\frac{1}{2\pi }\int_{-\infty
}^{\infty }\hat{f}(t)e^{iyt}dt)dy \\
&=&\frac{1}{2\pi }\int_{-\infty }^{\infty }\hat{f}(t)(\int_{-\infty
}^{\infty }g(y)e^{-ixy}e^{iyt}dy)dt \\
&=&\frac{1}{2\pi }\int_{-\infty }^{\infty }\hat{f}(t)\hat{g}(x-t)dt \\
&=&\frac{1}{2\pi }(\hat{f}\ast \hat{g})(x).
\end{eqnarray*}%
The change in the order of integration is justified by Fubini's theorem,
since due to boundedness of $\hat{f},$ $\hat{f}g\in L^{1}(\mathbb{R}).$ 
\qed 
\newline

We now investigate the question of uniqueness of Fourier transform, i.e, $%
\hat{f}=$ $\hat{g}$ implies $f=g.$ To show this, since Fourier transform is
a linear operator, it suffices to show that $\hat{f}=0$ implies $f=0$ a.e.

\begin{thm}[Uniqueness Theorem]
If $f\in L^{1}(\mathbb{R})$ and $\hat{f}=0$ everywhere ($\hat{f}$ is always
continuous), then $f=0$ a.e.
\end{thm}

\textbf{Proof:}\quad Let $e_{n}(x)$ be an approximate identity with compact
support and continuous derivative. By Theorem 1.6, $(e_{n}\ast f)\hat{\,}=%
\hat{e_{n}}\hat{f}=0$ everywhere. Since by Theorem 1.8, $e_{n}\ast f$ is
continuous and differentiable, by the inversion theorem, $e_{n}\ast f=0$
everywhere. But by Theorem 1.7, $e_{n}\ast f\rightarrow f$ in $L^{1}(\mathbb{%
R});$ so it follows that $f=0$ a.e. \qed \newline

\begin{df}
For $\mu \in M(\mathbb{R})$ (bounded Borel measure on $\mathbb{R}$, i.e., ${%
\displaystyle|\mu |(\mathbb{R})<\infty }$), define the Fourier-Stieltjes
transform $\mu \hat{\,}(y)$ as 
\begin{equation*}
\mu \hat{\,}(y)=\int_{-\infty }^{\infty }e^{ixy}d\mu (x).
\end{equation*}
\end{df}

Clearly, the Fourier-Stieltjes transform defines a bounded linear transform
from $M(\mathbb{R})$ to $\mathbb{C}.$

\begin{thm}[Uniqueness Theorem]
If $\hat{\mu}(y)=0$ for a.e. $y,$ then $\mu =0.$
\end{thm}

\textbf{Proof:}\quad Since $(C_{0}(\mathbb{R}))^{\ast }=M(\mathbb{R}),$ to
prove $\mu =0$ we need only to show that for all $h\in C_{0}(\mathbb{R}),$ ${%
\displaystyle\int h(t)d\mu (t)=0}.$ This is equivalent to showing that for
all $h\in C_{0}(\mathbb{R}),$ ${\displaystyle
(h\ast \mu )(0)=0},$ where ${\displaystyle(h\ast \mu )(x)=\int_{-\infty
}^{\infty }h(x-t)d\mu (t)}.$ Observe also that $h(x)\in C_{0}(\mathbb{R})$
if and only if $h(-x)\in C_{0}(\mathbb{R})$.

Assume that $\mu \hat{\,}=0.$ Then for all $f\in L^{1}(\mathbb{R}),$ $f\hat{%
\,}\ast \mu (x)=(f\ast \mu \hat{\,})(x)=0.$ Hence, if we prove that $\{f\hat{%
\,}:f\in L^{1}(\mathbb{R})\}$ is dense in $C_{0}(\mathbb{R})$, then for each 
$h\in C_{0}(\mathbb{R})$ there is $f_{n}\in L^{1}(\mathbb{R})$ such that $%
f_{n}\hat{\,}\rightarrow h$ in $C_{0}(\mathbb{R}).$ Since $f_{n}\hat{\,}\ast
\mu (x)\rightarrow h\ast \mu (x)$ at each $x,$ $h\ast \mu (x)=0.$

To show that $\{f\hat{\,}:f\in L^{1}(\mathbb{R})\}$ is dense in $C_{0}(%
\mathbb{R})$, we let 
\begin{equation*}
F(x)=\frac{1}{\sqrt{2\pi }}(\frac{\sin (x/2)}{x/2})^{2}
\end{equation*}%
and let $F_{\rho }(x)=\rho F(\rho x).$ Consider the integral 
\begin{equation*}
(h\ast F_{\rho })(x)=\frac{2}{\pi \rho }\int_{-\infty }^{\infty }h(x-u)\frac{%
\sin ^{2}(\rho u/2)}{u^{2}}du.
\end{equation*}%
Define 
\begin{equation*}
\mathcal{F}=\{(h\ast F_{\rho })(x):h\in C_{0}(\mathbb{R})\bigcap L^{1}(%
\mathbb{R});\,\rho >0\}.
\end{equation*}%
Clearly, $\mathcal{F}$ is a subset of $C_{0}(\mathbb{R})\bigcap L^{1}(%
\mathbb{R})$ and is dense in $C_{0}(\mathbb{R}).$

Let $h\in C_{0}(\mathbb{R})\bigcap L^{1}(\mathbb{R}).$ Then $(h\ast F_{\rho
})\hat{\,}(y)=h\hat{\,}(y)(F_{\rho })\hat{\,}(y).$ Since $h\in L^{1},$ $h%
\hat{\,}\in C_{0}(\mathbb{R}).$ Moreover, 
\begin{equation*}
F_{\rho }\hat{\,}(y)=\left\{ 
\begin{array}{ll}
1-\frac{|y|}{\rho } & \mbox{if $|y|\leq \rho$} \\ 
0 & \mbox{if $|y|>\rho$}%
\end{array}%
\right.
\end{equation*}%
belongs to $L^{1}(\mathbb{R}).$ Therefore, $(h\ast F_{\rho })\hat{\,}\in
L^{1}(\mathbb{R}).$ It follows from the inversion theorem that $h\ast
F_{\rho }$ is the Fourier transform of a function in $L^{1}(\mathbb{R}).$
Hence, $\mathcal{F}$ is a subset of $\{f\hat{\,}:f\in L^{1}(\mathbb{R})\}.$
Since $\mathcal{F}$ is dense in $C_{0}(\mathbb{R}),$ $\{f\hat{\,}:f\in L^{1}(%
\mathbb{R})\}$ is dense in $C_{0}(\mathbb{R}).$ \qed

\section{Kernels on $\mathbb{R}$}

We define the Dirichlet, Fej\'{e}r, and Poisson kernels on $\mathbb{R}$ by
defining their Fourier transforms, see H. Helson \cite{Helson}. 
\begin{equation*}
\hat{D}_{t}(y)=\left\{ 
\begin{array}{ll}
1 & \mbox{if $|y|\leq t$} \\ 
0 & \mbox{if $|y|>t$}%
\end{array}%
\right.
\end{equation*}%
\begin{equation*}
\hat{K}_{t}(y)=\left\{ 
\begin{array}{ll}
1-\frac{|y|}{t} & \mbox{if $|y|\leq t$} \\ 
0 & \mbox{if $|y|>t$}%
\end{array}%
\right.
\end{equation*}%
and 
\begin{equation*}
\hat{P}_{u}^{{}}(y)=e^{-u|y|}.
\end{equation*}%
The parameters $t$ and $u$ are positive, having limits $\infty $ and $0^{+},$
respectively.

Taking the inverse Fourier transform of $\hat{D}_{t}(y)$ we get the
Dirichlet kernel 
\begin{equation*}
D_{t}(x)=\frac{\sin tx}{\pi x}.
\end{equation*}%
Since $\hat{D}_{t}(y)\in L^{1}(\mathbb{R})$ and every point $y\neq t$ is a
Lipschitz point of $\hat{D}_{t}(y),$ it follows from the inversion theorem
that 
\begin{equation*}
\hat{D}_{t}(y)=\lim_{A,B\rightarrow \infty
}\int_{-B}^{A}D_{t}(x)e^{-ixy}dx,\,\,y\neq \pm t.
\end{equation*}%
That is, although $D_{t}(x)$ is not integrable, its Fourier transform in the
generalized sense is $\hat{D}_{t}(y).$ Since $\hat{D}_{t}(y)$ is
discontinuous, $D_{t}(x)$ cannot be integrable. Clearly, the Dirichlet
kernel does not belong to the family of approximate identities.

To calculate the Fej\'{e}r kernel, it follows from definitions that

\begin{eqnarray*}
(\hat{D}_{t}\ast \hat{D}_{t})(y) &=&\int_{-\infty }^{\infty }\hat{D}%
_{t}(y-\tau )\hat{D}_{t}(\tau )d\tau \\
&=&\int_{-t}^{t}\hat{D}_{t}(y-\tau )d\tau \\
&=&\int_{y-t}^{y+t}\hat{D}_{t}(u)du.
\end{eqnarray*}

To calculate the last integral, we consider two cases. If $|y|\geq 2t$, then
the intervals $[y-t,y+t]$ and $[-t,t]$ are disjoint so that the integral
equals zero; if $|y|<2t,$ then either $y+t$ or $y-t$ is in $(-t,t)$, but not
both, so that the integral equals $2t-|y|.$ Combining both results we get,

\vspace{0.1in}

$(\hat{D}_{t}\ast \hat{D}_{t})(y)=$%
\[
\int_{y-t}^{y+t}\hat{D}_{t}(u)du
=\left\{\begin{array}{ll}
2t-|y| &\mbox{if $|y|\leq 2t$}\\
$0$ &\mbox{if $|y|>2t$}
\end{array}
\right.
=2t\hat{K}_{2t}(y).
\]%

Also, by Theorem 1.14 we have that, 
\begin{eqnarray*}
(D_{t}\cdot D_{t})\hat{\,}(x) &=&\frac{1}{2\pi }(\hat{D}_{t}\ast \hat{D}%
_{t})(x) \\
&=&\frac{1}{2\pi }(2t\hat{K}_{2t}(x)).
\end{eqnarray*}

\bigskip

Therefore, it follows from the inversion theorem that $(D_{t}\cdot D_{t})(x)=%
\frac{1}{2\pi }(2tK_{2t}(x))$, or 
\begin{equation*}
2tK_{2t}(x)=\frac{1}{2\pi }(2\pi D_{t}(x))^{2}.
\end{equation*}

Hence, we obtain the Fej\'{e}r kernel%
\begin{equation*}
K_{t}(x)=\frac{1}{2\pi t}(\frac{\sin (\frac{tx}{2})}{\frac{x}{2}})^{2}.
\end{equation*}%
$K_{t}(x)$ is positive and integrable. Its Fourier transform is the function 
$\hat{K}_{t}(y)$ by the inversion theorem. Moreover, ${\displaystyle\int
K_{t}(x)dx=1}$ because $\hat{K}_{t}(y)=1$ at $y=0.$ For any $\epsilon >0,$ 
\begin{equation*}
\int_{|x|>\epsilon }K_{t}(x)dx\leq \frac{1}{2\pi t}\int_{|x|>\epsilon }\frac{%
4}{x^{2}}dx\rightarrow 0
\end{equation*}%
as $t\rightarrow \infty .$ Hence $(K_{t})$ is an approximate identity on $%
\mathbb{R}.$

A direct computation of the inverse Fourier transform of $\hat{P}_{u}(y)$
gives \quad 
\begin{eqnarray*}
P_{u}(x) &=&\frac{1}{2\pi }\int_{-\infty }^{\infty }\hat{P}_{u}(y)e^{ixy}dy
\\
&=&\frac{1}{2\pi }\int_{-\infty }^{\infty }e^{-u\left\vert y\right\vert
}e^{ixy}dy \\
&=&\frac{1}{2\pi }\int_{-\infty }^{0}e^{uy}e^{ixy}dy+\frac{1}{2\pi }%
\int_{0}^{\infty }e^{-uy}e^{ixy}dy \\
&=&\frac{1}{2\pi }\int_{-\infty }^{0}e^{y(u+ix)}dy+\frac{1}{2\pi }%
\int_{0}^{\infty }e^{-y(u-ix)}dy \\
&=&\frac{1}{2\pi }(\frac{1}{u+ix}+\frac{1}{u-ix}) \\
&=&\frac{u}{\pi (u^{2}+x^{2})}.
\end{eqnarray*}

This gives the formula for Poisson kernel

\begin{equation*}
P_{u}(x)=\frac{u}{\pi (u^{2}+x^{2})}.
\end{equation*}%
Clearly, $P_{u}$ is positive, and we check that 
\begin{equation*}
\lim_{u\rightarrow 0^{+}}\int_{-\epsilon }^{\epsilon }P_{u}(x)dx=1
\end{equation*}%
for each $\epsilon >0.$ Thus $(P_{u})$ is an approximate identity with $%
u\downarrow 0.$

\begin{thm}[Inversion Theorem]
If $f$ and $\hat{f}$ are both integrable, then $f(x)$ a.e. equals to a
continuous function which is the inverse Fourier transform of $\hat{f}$,
that is, 
\begin{equation}
f(x)=\frac{1}{2\pi }\int \hat{f}(y)e^{ixy}dy,\; a.e.
\end{equation}
\end{thm}

\textbf{Proof:}\quad $P_{u}\ast f$ is continuous. We have, 
\begin{eqnarray*}
(P_{u}\ast f)(x) &=&\int_{-\infty }^{\infty }P_{u}(x-t)f(t)dt \\
&=&\frac{1}{2\pi }\int_{-\infty }^{\infty }(\int_{-\infty }^{\infty
}e^{-u|y|}e^{iy(x-t)}dy)f(t)dt \\
&=&\frac{1}{2\pi }\int_{-\infty }^{\infty }e^{-u|y|}e^{ixy}(\int_{-\infty
}^{\infty }f(t)e^{-iyt}dt)dy \\
&=&\frac{1}{2\pi }\int_{-\infty }^{\infty }e^{-u|y|}e^{ixy}\hat{f}(y)dy.
\end{eqnarray*}%
Note that the $P_{u}\ast f$ converges to $f(x)$ in $L^{1}(\mathbb{R})$ so
that $P_{u}\ast f$ converges to $f$ almost everywhere at least on a
subsequence of $u\downarrow 0.$ We then obtain 
\begin{equation*}
f(x)=\lim_{u\downarrow 0}\frac{1}{2\pi }\int_{-\infty }^{\infty
}e^{-u|y|}e^{ixy}\hat{f}(y)dy=\frac{1}{2\pi }\int \hat{f}(y)e^{ixy}dy.
\end{equation*}%
The last limit holds because of Lebesgue's dominated convergence theorem. 
\qed

\begin{df}
For any $f\in L^p(\mathbb{R}), 1\leq p<\infty,$ we define the Poisson
integral of $f$ as 
\begin{equation}
F(x+iu)=P_u*f(x)=\frac{1}{\pi}\int_{\mathbb{R}} \frac{uf(s)}{u^2+(x-s)^2}ds.
\end{equation}
\end{df}

Since $P_{u}\in L^{q}(\mathbb{R}),$ $q$ conjugate exponent of $p,$ $F(x+iu)$
is defined as a continuous function of $x.$ \footnote{%
If $f\in L^{p}(\mathbb{R}),$ $1\leq p\leq \infty ,$ and $g\in L^{q}(\mathbb{R%
}),$ $\frac{1}{p}+\frac{1}{q}=1,$ then $(f\ast g)(x)$ exists everywhere,
belongs to $C(\mathbb{R}),$ and $||f\ast g||_{c}\leq ||f||_{p}||g||_{q}.$}
Moreover, $F(x+iu)$ provides a harmonic extension of $f$ to the upper half
plane. This can be verified directly.

\begin{thm}
The Poisson integral has a semigroup property: $P_u*P_v=P_{u+v}$ for all
positive $u$ and $v.$
\end{thm}

\textbf{Proof:}\quad We have that 
\begin{eqnarray*}
\hat{(P_u*P_v)}(y)&=&\hat{P_u}(y)\cdot\hat{P_v}(y)= e^{-u|y|}\cdot e^{-v|y|}
\\
&=&e^{-(u+v)|y|}=\hat{P_{u+v}}(y).
\end{eqnarray*}
It follows from the inversion theorem that $P_u*P_v=P_{u+v}.$ \qed

\begin{thm}
$||F(\cdot+iu)||_p$ increases as $u\downarrow 0,$ for any $p, 1\leq
p<\infty. $ (if $p=1$, consider $P_u*\mu$). Similarly, if $f$ is bounded, ${%
\displaystyle \sup_{x\in \mathbb{R}}|F(x+iu)|}$ increases as $u\downarrow 0.$
\end{thm}

\textbf{Proof:}\quad Let $v<u$ be given. Let $r=u-v\geq 0.$ Then 
\begin{equation*}
||P_{u}*f||_p=||P_{v+r}*f||_p= ||(P_r*P_v)*f||_p\leq ||P_r||_1||P_v*f||_p
=||P_v*f||_p. \qed
\end{equation*}


\begin{lem}
Let $f_u(x)=F(x+iu)$ be a harmonic function in the upper half plane such
that 
\begin{equation*}
\sup_{u>0}||f_u(\cdot)||_p=A<\infty.
\end{equation*}
Then 
\begin{equation*}
f_{u+v}(x)=(P_u*f_v)(x).
\end{equation*}
\end{lem}

\textbf{Proof:}\quad $f_{u+v}(x)=(P_{u}\ast f_{v})(x)$ says that the values
of $F(u+ix)$ at the level $u+v$ are the values of $F(u+ix)$ at the level $v$
convolved with the Poisson kernel with parameter $u.$ \footnote{%
In periodic case, $P_{r}\ast f_{s}=f_{rs}$ is proved by using the fact that
a harmonic function is the real part of an analytic function.}

We may assume that $F$ is real. Fix $v>0.$ Define $G(x+iu)=P_{u}\ast
f_{v}(x) $ ($G$ is the Poisson integral of the values of $F$ at level $v$). $%
G(x+iu)$ is harmonic in $u>0$ and $\sup_{u>0}||G(\cdot ,u)||_{p}\leq
||f_{v}(\cdot )||_{p}<\infty .$ Note that $G(x+iu)$ has boundary value
(pointwise limit) $f_{v}(x)$ as $u\rightarrow 0,$ which can be simply viewed
as the value of $G(x+iu)$ when $u=0.$ Therefore, $G(x+iu)-F(x+iu+iv)$ is a
harmonic function in $u>0,$ satisfying $\sup_{u>0}||G(\cdot +iu)-F(\cdot
+iu+iv)||_{p}<\infty , $ continuous on the closed upper half plane and null
on the real axis $u=0.$ Now, let 
\begin{equation*}
H(x+iu)=G(x+iu)-F(x+iu+iv).
\end{equation*}%
We must show that $H(x+iu)$ vanishes for $u>0.$

Let $h\in L^{1}(\mathbb{R})\bigcap L^{q}(\mathbb{R}),$ where $q$ is the
conjugate exponent of $p.$ Define 
\begin{equation*}
L(x+iu)=\int_{-\infty }^{\infty }h(x-y)H(y+iu)dy.
\end{equation*}%
Then $L(x+iu)$ is continuous on the closed upper half plane, harmonic in the
upper half plane, and is null on the real axis. Moreover, $%
\sup_{u>0}|L(x+iu)|\leq ||h||_{q}\sup_{u>0}||H(\cdot +iu)||_{p}<\infty ,$
that is, $L(x+iu)$ is bounded on the upper half plane. Extend this function
to a bounded harmonic function on the whole plane by setting $%
L(x-iu)=-L(x+iu)$ (Reflection Principle for Harmonic Functions). By
Liouville's theorem, $L$ is constant. Letting $h$ range over an approximate
identity shows that $H$ is a constant, and since it vanishes on the real
axis, is null. \qed

\begin{thm}
Let $f_u(x)=F(x+iu)$ be a harmonic function in the upper half plane. Then
there exists a $f\in L^p(\mathbb{R}),$ $1<p\leq \infty,$ so that $%
f_u(x)=P_u*f(x)$ if and only if $f_u(x)\in L^p(\mathbb{R})$ with the norm
bounded by a constant independent of $u>0,$ that is, 
\begin{equation*}
\sup_{u>0}A_u=\sup_{u>0}||f_u(\cdot)||_p=A<\infty.
\end{equation*}
\end{thm}

\textbf{Proof:}\quad \emph{Necessity:} \quad If we think of $P_{u}\ast f(x)$
as a family (with continuous parameter $u>0$) of functions $f_{u}(x)$
defined on $\mathbb{R},$ then as $p>1,$ 
\begin{equation*}
||f_{u}||_{p}=||P_{u}\ast f||_{p}\leq ||P_{u}||_{1}||f||_{p}.
\end{equation*}%
Hence, $\{f_{u}\},\; u >0 $, is bounded in $L^{p}(\mathbb{R}).$\newline

\emph{Sufficiency:}\quad Assume that $f_{u}(x)=F(x+iu)$ is bounded in $L^{p}(%
\mathbb{R}).$ If $1<p\leq \infty ,$ by Banach-Alaoglu's theorem (view $L^{p}(%
\mathbb{R})$ as the dual of separated normed space $L^{q}(\mathbb{R}),$ $%
1\leq q<\infty $), $\{f_{u}\} $ is weakly* sequentially compact in $L^{p}(%
\mathbb{R}),$ that is, there is an element $f $ of $L^{p}(\mathbb{R})$ such
that every *-neighborhood of $f$ contains $f_{u}$ for arbitrary small
positive $u.$ In other words, there is a subsequence $f_{v_{j}}$ of $f_{u}$
that is weakly* convergent to \emph{some} $f\in L^{p}(\mathbb{R})$ as $%
v_{j}\downarrow 0^{+},$ i.e., for all $g\in L^{q}(\mathbb{R}), $ $%
\displaystyle{\int f_{v_{j}}g\rightarrow \int fg}$ as $v_{j}\downarrow
0^{+}. $ In particular, since for each $x,$ $P_{u}(x-t)\in L^{q}(\mathbb{R}%
), $ $1\leq q<\infty ,$ we have $P_{u}\ast f_{v_{j}}(x)=\displaystyle{\int
P_{u}(x-t)f_{v_{j}}(t)dt}$ tends to $\displaystyle{\int
P_{u}(x-t)f(t)dt=P_{u}\ast f(x)}$ as $v_{j}\downarrow 0^{+}.$ On the other
hand, $P_{u}\ast f_{v_{j}}(x)=f_{u+v_{j}}(x)$ (see Lemma 2.1), which
converges to $f_{u}(x)$ by the continuity of $F(x+iu).$ Therefore, $%
P_{u}\ast f(x)=f_{u}(x)$ for all $x.$

If $1<p<\infty ,$ then $P_{u}\ast f\rightarrow f$ in the norm of $L^{p}(%
\mathbb{R})$ (Fejer's theorem). If $p=\infty $ then $P_{u}\ast f\rightarrow
f $ in weak* topology in $L^{\infty }(\mathbb{R}),$ i.e., for every $s(x)\in
L^{1}(\mathbb{R}),$ $\lim_{u\downarrow 0^{+}}\int [P_{u}\ast
f(x)-f(x)]s(x)dx=0.$ (For a proof, see Butzer \cite{Butzer&Nessel}.) \qed 
\newline

\begin{thm}
Let $f_u(x)=F(x+iu)$ be a function harmonic in the upper plane $u>0.$ Then
there is a unique measure $\mu\in M(\mathbb{R})$ such that 
\begin{equation*}
f_u(x)=F(x+iu)=P_r*\mu(x)=\int_{-\infty}^\infty P_u(x-t)d\mu(t)
\end{equation*}
if and only if 
\begin{equation*}
A_u=\int |F(x+iu)|dx\leq K,\quad \forall\quad u>0.
\end{equation*}
Moreover, $||\mu||=\lim_{u\downarrow 0}A_u.$
\end{thm}

\textbf{Proof:}\quad \emph{Necessity:} \quad If we think of $P_u*\mu(x)$ as
a family (with continuous parameter $u>0$) of functions defined on $\mathbb{R%
}$, then 
\begin{equation*}
||f_u||_1=||P_u*\mu||\leq ||P_u||_1||\mu||.
\end{equation*}
Therefore, $\{f_u\}, u>0,$ is bounded in $L^1(\mathbb{R}).$\newline

\emph{Sufficiency:}\quad By assumption, $||f_u||_1\leq K,$ i.e., $%
||f_u(x)dx||_{M(\mathbb{R})}=||f_u||_1\leq K, \quad \forall u>0.$ Since $C_0(%
\mathbb{R}),$ as the pre-dual of $M(\mathbb{R}),$ is separable normed space,
by Banach-Alaoglu theorem the closure of $\{f_u(x)dx\}$ in $M(\mathbb{R})$
is weak* sequentially compact. Therefore, there is a subsequence $%
\{f_{v_j}\}(x)dx$ of $f_u(x)dx$ that converges to some $\mu\in M(\mathbb{R})$
in weak* topology. That is, 
\begin{equation*}
\int h(e^{-it})f_{v_j}(t)dt\rightarrow \int h(e^{-it})d\mu(t),\qquad
v_j\rightarrow 0
\end{equation*}
for each $h\in C_0(\mathbb{R}).$ In particular, since for each $x,$ $%
P_u(x-t)\in C_0(\mathbb{R}),$ 
\begin{equation*}
\int P_u(x-t)f_{v_j}(t)dt\rightarrow \int P_u(x-t)d\mu(t), \,\,
v_j\rightarrow 0.
\end{equation*}
On the other hand, 
\begin{equation*}
P_u*f_{v_j}(x)=f_{u+v_j}(x)\rightarrow f_u(x),\,\, v_j\rightarrow 0.
\end{equation*}
Hence, $f_u(x)=\displaystyle{\int P_u(x-t)d\mu(t)}$ for all $x.$

We show that $||\mu||=\lim_{u\downarrow 0}A_u.$ Note that $%
\mu=\lim_{j\rightarrow \infty}f_{v_j}(x)dx$ in the weak* topology of $M(%
\mathbb{R})$ as the dual of $C_0(\mathbb{R}).$ It follows that $||\mu||\leq
\liminf_{j\rightarrow \infty} A_{v_j}$ where $A_{v_j}=||f_{v_j}||_1$ (For a
proof, see the Appendix). Since $A_u$ increases with $u\downarrow $ and $%
A_u\leq K,$ $||\mu||\leq \lim_{u\rightarrow 0}A_u.$ Furthermore, the
inequality cannot be strict. Note that $f_u=P_u*\mu$ and $||f_u||_1\leq
||P_u||_1||\mu||.$ Therefore, $A_u=||f_u||_1\leq ||\mu||$ for every $u>0.$
If the inequality were strict, we would have $A_u\leq
||\mu||<\lim_{u\rightarrow 0}A_u$ for $u>0,$ which is impossible.

As to the norm convergence of $||f_{u}-\mu ||_{M(\mathbb{R})}\rightarrow 0$
as $u\rightarrow 0,$ if $\mu $ is absolutely continuous then $\mu =f(x)dx$
for some $f\in L^{1}(\mathbb{R}).$ Hence $f_{u}=P_{u}\ast \mu $ is indeed $%
f_{u}=P_{u}\ast f.$ Thus, by Fejer's theorem, $||f_{u}-f||_{1}\rightarrow 0.$
That is, $||f_{u}-\mu ||_{M(\mathbb{R})}\rightarrow 0$ as $u\rightarrow 0.$ 
\qed

\section{The Plancherel Theorem}

In this section we define 
\begin{equation*}
\hat{f}(y)=\frac{1}{\sqrt{2\pi}}\int_{-\infty}^\infty f(x)e^{-ixy} dx.
\end{equation*}

\vspace{0.1in}

\begin{lem}
Let $\mathcal{C}$ be the collection of continuously differentiable functions
with compact support. Then $\mathcal{C}\subset L^1(\mathbb{R})\bigcap L^2(%
\mathbb{R})$ and $\mathcal{C}$ is a dense subspace of $L^2(\mathbb{R}).$
\end{lem}

\textbf{Proof:}\quad Let $f\in L^2(\mathbb{R}).$ Define $f_k(x)=f(x)$ if $%
|x|\leq k$; and $f_k(x)=0$ if $|x|>0.$ Then $f_k\rightarrow f$ in $L^2(%
\mathbb{R}).$ Furthermore, we may choose an approximate identity with
compact support and continuous derivative, for instance, let $h(x)=%
\displaystyle{e^{-\frac{1}{x^2}}}$ for $x\geq 0$ and $h(x)=0$ for $x<0.$
Then $h\in C^\infty(\mathbb{R})$ and $\phi(x)=h(x+1)h(1-x)\in C^\infty(%
\mathbb{R})$ and has compact support $[-1,1],$ and $\int \phi(x)dx=1,$ when
properly normalized. Let $e_n(x)=n\phi(nx).$ Then $e_n(x)$ is an approximate
identity with compact support and continuous derivative (in fact, $C^\infty(%
\mathbb{R})$). Since for each $k,$\, $f_k$ has compact support, $e_n*f_k$
provides a continuously differentiable approximation with compact support to 
$f_k$ in $L^2(\mathbb{R}).$ Hence $\mathcal{C}$ is dense in $L^2(\mathbb{R}%
). $

\begin{lem}
If $f\in \mathcal{C}$, then $\hat{f}\in L^{2}(\mathbb{R}).$ Moreover, $||%
\hat{f}||_{2}=||f||_{2}.$ Hence, the Fourier transform $\hat{f}$ (as defined
in this section) is isometric from $\mathcal{C}$ to $\mathcal{F}(\mathcal{C}%
) $ as subspaces of $L^{2}(\mathbb{R}).$
\end{lem}

\textbf{Proof:}\quad Let $f\in \mathcal{C}.$ Define $\tilde{f}(x)=\overline{%
f(-x)}. $ Then $f\ast \tilde{f}(x)\in \mathcal{C}.$ By the inversion
theorem, at every point $x$ where $(f\ast \tilde{f})(x)$ satisfies the
Lipschitz condition, we have 
\begin{equation*}
(f\ast \tilde{f})(x)=\lim_{A,B\rightarrow \infty }\frac{1}{\sqrt{2\pi }}%
\int_{-B}^{A}\widehat{f\ast \tilde{f}}(y)e^{ixy}dy.
\end{equation*}%
Since $f\ast \tilde{f}(x)\in \mathcal{C},$ it satisfies the Lipschitz
condition at \emph{every} point, in particular, at $x=0,$ we have 
\begin{equation*}
(f\ast \tilde{f})(0)=\lim_{A,B\rightarrow \infty }\frac{1}{\sqrt{2\pi }}%
\int_{-B}^{A}\widehat{f\ast \tilde{f}}(y)dy.
\end{equation*}%
Note that $||f||_{2}=(f\ast \tilde{f})(0),$ $\widehat{\tilde{f}}=\overline{%
\hat{f}},$ and $\displaystyle{\frac{1}{\sqrt{2\pi }}\widehat{f\ast \tilde{f}}%
(y)=|\hat{f}(y)|^{2}}.$ We have $||f||_{2}=||\hat{f}||_{2}.$\newline

The Fourier transform is an isometry defined on $\mathcal{C}.$ Since it is
defined on a dense subspace of $L^2(\mathbb{R}),$ it has a unique \emph{%
continuous} extension to an isometry $\mathcal{F}$ of all of $L^2(\mathbb{R}%
) $ \emph{into} itself, which is defined as follows: for $f\in L^2(\mathbb{R}%
), $ let $f_n\in \mathcal{C}$ such that $f_n\rightarrow f.$ Since $||%
\mathcal{F}f||_2=||f||_2$ for all $f\in \mathcal{C},$ $\hat{f_n}$ is a
Cauchy sequence in $L^2(\mathbb{R})$ and so converges to \emph{some} $g\in
L^2(\mathbb{R}).$ We define $\mathcal{F}(f)=g.$ Let us show $||\mathcal{F}%
(f)||_2=||f||_2$ for all $f\in L^2(\mathbb{R}).$ Let $f\in L^2(\mathbb{R})$
and $f_k\in \mathcal{C}\rightarrow f.$ Then by definition of $\mathcal{F},$ $%
||\mathcal{F}f_k||_2\rightarrow ||\mathcal{F}f||.$ On the other hand, $||%
\mathcal{F}f_k||=||f_k||_2\rightarrow ||f||_2.$ Therefore, $\mathcal{F}$ is
an isometry of $L^2(\mathbb{R})$ \emph{into} $L^2(\mathbb{R}).$ We will
prove that $\mathcal{F}$ is indeed `onto'. \qed

\begin{lem}
The Fourier transform of $L^2(\mathbb{R})$ is onto, i.e., $E=\{\mathcal{F}%
(f): f\in L^2(\mathbb{R})\}=L^2(\mathbb{R}).$
\end{lem}

\textbf{Proof:}\quad First we prove that $E$ is dense in $L^2(\mathbb{R}).$

We prove that for each $h\in \mathcal{C},$ $<\mathcal{F}f,h>=<f,h^*>$ for
all $f\in L^2(\mathbb{R}),$ where $h^*$ is defined by the formula: 
\begin{equation*}
h^*(x)=\frac{1}{\sqrt{2\pi}}\int h(y)e^{ixy}dy.
\end{equation*}
In fact, for $f,h\in \mathcal{C}$ we have 
\begin{equation*}
\int \hat{f}\overline{h}=\int f\overline{h^*},
\end{equation*}
that is, $<\mathcal{F}f,h>=<f,h^*>$ for all $f\in \mathcal{C}.$ It follows
that $<\mathcal{F}f,h>=<f,h^*>$ for all $f\in L^2(\mathbb{R}).$

The operator $\mathcal{F}^*$ defined on $L^2(\mathbb{R})$ by $\mathcal{F}%
^*h=h^*$ is called the \emph{adjoint operator} of $\mathcal{F}.$ (see the
Appendix). Note that $\mathcal{F} $ is essentially the Fourier transform,
and therefore, is an isometry. Thus its null space $N(\mathcal{F}^*)$
contains $0$ only (uniqueness theorem for F.T.). Since $N(T^*)^{\perp}=R(T)$
(see the Appendix), it follows that $E,$ the range of $\mathcal{F},$ is
dense in $L^2(\mathbb{R})$.

To prove $E=L^2(\mathbb{R}),$ we show that $E$ is closed. Take $g\in 
\overline{E}.$ Then there exists $g_k\in E$ with $g_k\rightarrow g.$ Let $%
f_k $ be such that $\mathcal{F}f_k=g_k.$ Since $\mathcal{F}$ is an isometry, 
$f_k $ is a Cauchy sequence converging to some $f\in L^2(\mathbb{R})$ and we
must have $\mathcal{F}f=g.$ Since $\overline{E}=L^2(\mathbb{R})$ and $E$ is
closed, $E=L^2(\mathbb{R}).$ \qed \newline

\begin{thm}[Plancherel]
The Fourier transform $\mathcal{F}$ is a unitary operator on $L^2(\mathbb{R}%
) $ and the inverse Fourier transform, $\mathcal{F}^{-1},$ can be obtained
by $(\mathcal{F}^{-1}f)(x)=(\mathcal{F}f)(-x)$ for all $f\in L^2(\mathbb{R}%
). $
\end{thm}

\textbf{Proof:}\quad Since $\mathcal{F}$ is an isometry of $L^{2}(\mathbb{R}%
) $ onto $L^{2}(\mathbb{R}),$ $\mathcal{F}$ is a unitary operator on $L^{2}(%
\mathbb{R}).$ It follows from the properties of a unitary operator that $%
\mathcal{F}^{-1}=\mathcal{F}^{\ast }$ (see the Appendix). The form of $%
\mathcal{F}^{\ast }$ can be easily found when acting on $f\in \mathcal{C}:$
Let $g\in \mathcal{C}.$ A change of order of integration gives 
\begin{equation*}
\displaystyle{\int (\mathcal{F}g)(x)\overline{f(x)}dx=\int g(x)\overline{%
f^{\ast }(x)}dx,}
\end{equation*}%
i.e., $<\mathcal{F}(g),f>=<g,f^{\ast }>$ for all $g\in \mathcal{C},$ where $%
f^{\ast }=\mathcal{F}^{\ast }(f)$ is defined by the formula: 
\begin{equation*}
f^{\ast }(x)=\frac{1}{\sqrt{2\pi }}\int_{R}f(y)e^{ixy}dy.
\end{equation*}%
It follows that $(\mathcal{F}^{-1}f)(x)=(\mathcal{F}f)(-x)$ for all $f\in 
\mathcal{C}. $ For $f\in L^{2}(\mathbb{R}),$ we take $f_{k}\in \mathcal{C}$
with $f_{k}\rightarrow f$ in $L^{2}(\mathbb{R}).$ Then 
\begin{equation*}
(\mathcal{F}^{-1}f)(x)=l.i.m.(\mathcal{F}^{-1}f_{k})(x)=l.i.m.(\mathcal{F}%
f_{k})(-x)=(\mathcal{F}f)(-x).
\end{equation*}%
This shows that $(\mathcal{F}^{-1}f)(x)=(\mathcal{F}f)(-x)$ for all $f\in
L^{2}(\mathbb{R}).$ \qed \newline

\begin{lem}[Multiplication Formula]
If $f,g\in L^1(\mathbb{R}),$ then 
\begin{equation}
{\displaystyle\int \hat{f}g=\int f\hat{g}}.
\end{equation}
\end{lem}

\textbf{Proof:}\quad Since $F(x,t)=f(t)g(x)e^{-itx}$ is a measurable
function on $\mathbb{R}\times \mathbb{R}$ and $||F(x,t)||_{L(\mathbb{R}%
^{2})}=||f||_{1}||g||_{1},$ we can apply Fubini's theorem to obtain 
\begin{equation*}
\int \hat{f}(t)g(t)dt=\int f(t)(\int g(x)e^{-ixt}dx)dt=\int f(t)\hat{g}%
(t)dt. \qed
\end{equation*}

As an application of the multiplication formula, we prove the following
Fourier inverse theorem.

\begin{lem}
If $f\in L^1(\mathbb{R}),$ then the Abel mean of the Fourier integral
converges to $f(x)$ a.e., i.e., for almost every $x,$ 
\begin{equation}  \label{eq}
\lim_{\epsilon\rightarrow 0}\frac{1}{2\pi} \int \hat{f}(t)e^{-\epsilon
|t|}e^{ixt}dt=f(x).
\end{equation}
\end{lem}

\textbf{Proof:}\quad Let 
\begin{equation*}
g(t)=e^{-\epsilon |t|}e^{ixt},\,\, \epsilon>0.
\end{equation*}
Then 
\begin{equation*}
\hat{g}(t)=2\pi P_\epsilon(x-t)=\frac{2\epsilon}{\epsilon^2+ (x-t)^2}.
\end{equation*}
Using the multiplication formula, we get, if $f\in L^1(\mathbb{R}),$ 
\begin{equation*}
\frac{1}{2\pi} \int \hat{f}(t)e^{-\epsilon |t|}e^{ixt}dt=\int
f(t)P_\epsilon(x-t)dt.
\end{equation*}
Since the latter convolution converges to $f(x)$ a.e., \footnote{%
At every point of $x$ for which 
\begin{equation*}
\int_0^h[f(x+u)+f(x-u)-2f(u)]du={\scriptsize {o}(h),}
\end{equation*}
$\lim f(t)P_\epsilon(x-t)dt\rightarrow f(x).$} 
\begin{equation*}
\lim_{\epsilon\rightarrow 0} \frac{1}{2\pi}\int \hat{f}(t)e^{-\epsilon
|t|}e^{ixt}dt=f(x). \qed
\end{equation*}
\newline

\begin{cor}[Inversion Theorem]
If $f\in L^{1}(\mathbb{R})$ so that $\hat{f}\in L^{1}(\mathbb{R}),$ then for
a.e. $x,$ 
\begin{equation}
\frac{1}{2\pi }\int \hat{f}(t)e^{ixt}dt=f(x).
\end{equation}
In particular, the inversion formula holds at every $x$ for which 
\begin{equation*}
\int_{0}^{h}[f(x+u)+f(x-u)-2f(u)]du=o(h)
\end{equation*}%
holds.
\end{cor}

\textbf{Proof:}\quad If $\hat{f}\in L^1(\mathbb{R})$, the corollary follows
from Lemma 3.5 by applying the Lebesgue Dominated Convergence Theorem. 
\qed 
\newline

\begin{lem}
If $f(x)\in L^1(\mathbb{R})$ is continuous at $x=0$ such that $\hat{f}\geq
0, $ then $\hat{f}\in L^1(\mathbb{R})$ and 
\begin{equation*}
f(x)=\frac{1}{2\pi}\int \hat{f}(y)e^{ixy}dy, a.e.
\end{equation*}
In particular, 
\begin{equation*}
f(0)=\frac{1}{2\pi}\int \hat{f}(y)dy.
\end{equation*}
\end{lem}

\textbf{Proof:}\quad We need only to show that $\hat{f}\in L^1(\mathbb{R}).$
Then the rest of the statements follows from the inversion theorem.

Note that Corollary 3.1 holds at every point $x$ for which 
\begin{equation*}
\int_0^h[f(x+u)+f(x-u)-2f(u)]du=o(h).
\end{equation*}
In particular, it holds at the point $x=0$ of continuity of $f,$ i.e., 
\begin{equation*}
\frac{1}{2\pi}\int \hat{f}(t)e^{-\epsilon |t|}dt=f(0).
\end{equation*}
By Fatou's lemma, 
\begin{equation*}
\int \hat{f}(y)dy\leq \lim_{\epsilon\rightarrow 0} \int \hat{f}%
(t)e^{-\epsilon |t|}dt=f(0).
\end{equation*}
Since $0\leq \hat{f},$ $\hat{f}\in L^2(\mathbb{R}).$ \qed \newline

\begin{lem}
If $f\in L^1(\mathbb{R})\bigcap L^2(\mathbb{R}),$ then $\hat{f}\in L^2(%
\mathbb{R})$ and $||\hat{f}||_2=||f||_2.$
\end{lem}

\textbf{Proof:}\quad Define $\tilde{f}(x)=\overline{f(-x)}.$ Since $f, 
\tilde{f}\in L^2(\mathbb{R}),$ $h=f*\tilde{f}\in L^1(\mathbb{R})$ and is
continuous. Further, $\hat{h}(y)=|\hat{f}(y)|^2\geq 0.$ Hence, $\hat{h}\in
L^1(\mathbb{R})$ and $h(0)=\int \hat{h}(y)dy.$ It follows that 
\begin{equation*}
\int |\hat{f}(y)|^2dy=\int \hat{h}(y)dy=h(0)= f*\tilde{f}(0)=\int
|f(x)|^2dx. \qed \newline
\end{equation*}

Since the Fourier transform is an isometry of $L^1(\mathbb{R})\bigcap L^2(%
\mathbb{R}),$ it has a unique \emph{continuous} extension to an isometry $%
\mathcal{F}$ of all of $L^2(\mathbb{R})$ \emph{into} itself with $||\mathcal{%
F}f||_2=||f||_2$ for all $f\in L^2(\mathbb{R}).$ For $f\in L^2(\mathbb{R}),$
define $f_n(x)=f(x)$ if $|x|\leq n$ and $f_n(x)=0$ if $|x|>n. $ Then $f_n\in
L^1(\mathbb{R})\bigcap L^2(\mathbb{R})$ and $f_n\rightarrow f$ in $L^2(%
\mathbb{R}).$ Define $\mathcal{F}f=l.i.m. \hat{f_n}.$ We'll prove that $%
\mathcal{F}$ is indeed `onto'.

\begin{lem}[Multiplication Formula for $L^2(\mathbb{R})$]
If $f, g\in L^2((\mathbb{R})$ then 
\begin{equation}
{\displaystyle \int \hat{f}g=\int f\hat{g}}.
\end{equation}
\end{lem}

\textbf{Proof:}\quad Fix $g\in L^{1}(\mathbb{R})\bigcap L^{2}(\mathbb{R})$
first. Let $f\in L^{2}(\mathbb{R})$ and $f_{k}\in L^{1}(\mathbb{R})\bigcap
L^{2}(\mathbb{R})$ with $l.i.m.f_{k}=f.$ Since $\hat{g}\in L^{2}(\mathbb{R}%
), $ ${\displaystyle\int f_{k}\hat{g}\rightarrow \int f\hat{g}}.$ It follows
from the multiplication formula for $L^{1}(\mathbb{R})$ that ${\displaystyle%
\int f_{k}\hat{g}=\int \hat{f_{k}}g\rightarrow \int \hat{f}g}.$ Hence for $%
f\in L^{2}(\mathbb{R})$ and $g\in L^{1}(\mathbb{R})\bigcap L^{2}(\mathbb{R}%
), $ ${\displaystyle\int \hat{f}g=\int f\hat{g}}.$ Starting with this
formula, for $f,g\in L^{2}(\mathbb{R}),$ we approximate $g$ by $g_{k}\in
L^{1}(\mathbb{R})\bigcap L^{2}(\mathbb{R}).$ \qed \newline

\begin{thm}[Plancherel]
The Fourier transform $\mathcal{F}$ is a unitary operator of $L^2(\mathbb{R}%
) $ and the inverse Fourier transform, $\mathcal{F}^{-1},$ can be obtained
by $(\mathcal{F}^{-1}f)(x)=(\mathcal{F}f)(-x)$ for all $f\in L^2(\mathbb{R}%
). $
\end{thm}

\textbf{Proof:}\quad We have already proved that $\mathcal{F}$ is an
isometry, we only need to show $\mathcal{F}$ maps $L^2(\mathbb{R})$ \emph{%
onto} $L^2(\mathbb{R}),$ i.e., $E=\{\mathcal{F}(f): f\in L^2(\mathbb{R}%
)\}=L^2(\mathbb{R}).$ As proven before, $E$ is closed. Assume that $E\neq
L^2(\mathbb{R}).$ Then there exists $g\neq 0, g\in L^2(\mathbb{R})\setminus
E,$ such that $<g, f>=0$ for all $f\in E,$ or $<g, \hat{h}>=0$ for all $h\in
L^2(\mathbb{R}). $ It follows from the multiplication formula that ${%
\displaystyle\int h \overline{\hat{g}}=0}$ for all $h\in L^2(\mathbb{R}).$
In particular, taking $h=\hat{g}\in L^2(\mathbb{R}),$ $||\hat{g}%
||_2=0=||g||_2$ and $g=0$ a.e., contrary to the assumption $g\neq 0.$
Therefore, $\mathcal{F}$ is onto and so is a unitary operator of $L^2(%
\mathbb{R}).$ \qed \newline

\section{Appendix}

\subsection{Weak/Weak * Topologies in Linear Spaces}

Let $X$ be a topological linear space and $X^{\prime}$ be its conjugate
space of all continuous linear functionals on $X.$ \footnote{%
When $X$ is a Hausdorff locally convex space, the Hahn-Banach theorem
ensures the existence of enough elements in $X^{\prime}$ to make possible a
rich theory of the duality between $X$ and $X^{\prime}.$}\newline

The \emph{weak topology} $\sigma(X, X^{\prime})$ on $X$ is defined as
follows:

Let $F$ be a nonempty finite subset of $X^{\prime}.$ Define 
\begin{equation*}
p_F(x)=max_{x^{\prime}\in F} |x^{\prime}(x)|,\quad x\in X.
\end{equation*}
$p_F(x)$ is a seminorm on $X.$\,\, $\sigma(X,X^{\prime})$ is the locally
convex topology on $X$ defined by the family of all seminorms $p_F(x),$
where $F$ ranges over all finite subsets of $X^{\prime}.$ A base at $x_0\in
X $ for this topology is given by sets of the form 
\begin{eqnarray*}
U_{F,r}&=&\{x :|x^{\prime}(x)-x^{\prime}(x_0)|<r \quad%
\mbox{for each $x'\in
F$}\} \\
&=&\bigcap_{x^{\prime}\in F}\{x: |x^{\prime}(x)-x^{\prime}(x_0)|<r\},
\end{eqnarray*}
where $r>0$ and $F$ is a nonempty finite subset of $X^{\prime}.$\,\, $%
\sigma(X,X^{\prime})$ is the weakest topology on $X$ for which all the
elements of $X^{\prime}$ are continuous.\newline

A sequence $\{x_{n}\}$ in a normed linear space $X$ converges to an element $%
f\in X^{\prime }$ in weak topology if and only if $\lim_{n\rightarrow \infty
}f(x_{n})=f(x)$ for all $f\in X^{\prime }.$\newline

The \emph{weak* topology} $\sigma(X^{\prime}, X)$ on $X^{\prime}$ is defined
as follows:

Let $A$ be a nonempty finite subset of $X.$ Define 
\begin{equation*}
p_A(x^{\prime})=max_{x\in A} |x^{\prime}(x)|,\quad x^{\prime}\in X^{\prime}.
\end{equation*}
$p_A(x^{\prime})$ is a seminorm on $X^{\prime}.$ $\sigma(X^{\prime},X)$ is
the locally convex topology on $X^{\prime}$ defined by the family of all
seminorms $p_A(x^{\prime}),$ where $A$ ranges over all finite subsets of $X.$
A base at $x_0^{\prime}\in X^{\prime}$ for this topology is given by sets of
the form 
\begin{eqnarray*}
U_{A,r}&=&\{x^{\prime}:|x^{\prime}(x)-x^{\prime}_0(x)|<r \quad%
\mbox{for each
$x\in A$}\} \\
&=&\bigcap_{x\in A}\{x: |x^{\prime}(x)-x^{\prime}_0(x)|<r\},
\end{eqnarray*}
where $r>0$ and $A$ is a nonempty finite subset of $X.$ $\sigma(X^{%
\prime},X) $ is the weakest topology on $X^{\prime}$ for which $%
x^{\prime}(x),$ as a linear functional acting on $X^{\prime}$, is continuous.%
\newline

A sequence $\{f_{n}\}$ in $X^{\prime }$ of a normed linear space $X$
converges to an element $f\in X^{\prime }$ in weak* topology if and only if
at each $x$, $\lim_{n\rightarrow \infty }f_{n}(x)=f(x),$ see K. Yosida \cite%
{Yosida}.

\begin{thm}
If $X$ is a Banach space, then $\{f_{n}\}\subset X^{\prime }$ converges
weakly* to an element $f\in X^{\prime }$ if and only if (1). $\{||f_{n}||\}$
is bounded; and (2). $\lim f_{n}(x)=f(x)$ for all $x$ in a dense subset
(with respect to norm topology) of $X.$
\end{thm}

\begin{prob}
In $X=l^2,$ for $1\leq m<n<\infty,$ let $x_{mn}\in l^2$ be defined as 
\begin{equation*}
x_{mn}^{(k)}=\left\{%
\begin{array}{ll}
1 & \mbox{if $k=m$} \\ 
m & \mbox{if $k=n$} \\ 
0 & \mbox{otherwise}%
\end{array}
\right.
\end{equation*}
and let $A=\{x_{mn}: 1\leq m<n<\infty\}.$ Then no sequence of elements of $A 
$ converges weakly to the origin, yet the origin is an accumulation point of 
$A$ in the weak topology.
\end{prob}

\textbf{Proof:}\quad Note that $l^{2}$ is reflexible, so the weak topology
and the weak* topology coincide on $l^{2}.$ To prove that $0\in l^{2}$ is an
accumulation point of $A$ in the weak topology, i.e., to prove that, since $%
0\not\in A,$ for any weak neighborhood $S$ of $0,$ $S\bigcap A$ is not
empty, \footnote{%
If $A\subset X,$ then $x_{0}\in X$ is called an \textbf{accumulation point}
of $A$ if every neighborhood of $x_{0}$ contains a point of $A\setminus
\{x_{0}\}.$ If $A$ is a subset of a Hausdorff space $X$ and $x_{0}$ is an
accumulation point of $A,$ then every neighborhood of $x_{0}$ contains 
\textbf{infinitely} many points of $A.$ The \textbf{closure} of $A$ consists
of points $x$ such that every neighborhood of $x$ contains at least a point
of $A.$} we have to be able to write down a base at $0:$ 
\begin{equation*}
U_{F}=\bigcap_{b\in F}\{a\in l^{2}:|<a,b>|=|\sum a^{(k)}\overline{b^{(k)}}%
|<\epsilon \},
\end{equation*}%
where $F$ is a finite subset of $l^{2}$ and $\epsilon >0.$ In particular,
for every fixed $b\in l^{2},$ 
\begin{equation*}
U_{b,\epsilon }(0)=\{a\in l^{2}:|<a,b>|<\epsilon \}
\end{equation*}%
is a weak neighborhood of $0.$

Given a weak neighborhood $U=U_{b,\epsilon }(0)$ of $0$, can we always find $%
x_{mn}\in A$ so that $x_{mn}\in U$? Observe that $%
|<x_{mn},b>|=|b^{(m)}+mb^{(n)}|,$ which can be made as small as we wish.
First we choose $m$ large enough so that $|b^{(m)}|$ is very small, then for
this fixed $m,$ choose $n$ large enough so that $|mb^{(n)}|$ is also very
small.\newline

Can we prove that there is no sequence of elements in $A$ that converges
weakly to $0$? Given any sequence of elements in $A,$ we show that there
exist $\epsilon_0>0$ and $b\in l^2$ (i.e. there exists a weak neighborhood $%
U_{b,\epsilon_0}(0)$ of $0$) such that for any $l,$ we can always find an
element $a$ in this sequence with subscript $\geq l$ such that $a\not\in
U_{b,\epsilon_0}(0).$

Consider a sequence, $\xi $, of elements of $x_{mn}\in A.$ If \emph{some
integer,} say $l,$ appears infinitely many times as the $m$-index of $%
x_{mn}\in \xi ,$ then we choose $b$ so that $b^{(l)}=1,$ $b^{(k)}=0$ $k\neq
l.$ Of course, $b\in l^{2}$ and there is a (of course, infinite) subsequence 
$\{x_{ln}\}$ of $\xi $ with $|<b,x_{ln}>|=1.$ If none of the integers
appears infinitely many times as $m$-index in $\xi ,$ then the range of $m$-
index of elements $x_{mn}\in \xi $ is unbounded. We may extract a
subsequence, call it $\eta ,$ of $\xi $ so that their $m$-indices form a
(strictly) increasing sequence. Note that the range of $n$-index of $%
x_{mn}\in \eta $ is unbounded. We may extract a further subsequence, call it 
$\zeta ,$ of $\eta $ so that their $n$-indices form a (strictly) increasing
sequence of integers. Now we define $b$ with $b^{(n)}=\frac{1}{m}$ if $%
x_{mn}\in \zeta $ (Note: for each $n$ there is only one $m$ such that $%
x_{mn}\in \zeta $) and $b^{(n)}=0,$ otherwise. Note that if $x_{mn}\in \zeta
,$ then $|<x_{mn},b>|=1. $ All that remains is to notice that $b\in l^{2}.$

\begin{df}
Let $X$ be a topological space. If $A\subset X$ is such that every sequence
in $A$ has a subsequence that converges to a point in $A$, then $A$ is
called \textbf{sequentially compact.}
\end{df}

\begin{thm}
Let $X$ be a normed linear space. If $F\subset X^{\prime}$ is weak*
sequentially compact, then $F$ is countably weak* compact.
\end{thm}

\textbf{Proof:}\quad Suppose that there is an open cover (in weak* topology) 
$\{U_j\}$ of $F$ for which there is no finite subcover. Then for any finite
collection $U_j,$ $1\leq j\leq n,$ $F\setminus \bigcup_{j=1}^nU_j\neq \phi.$
Pick $x_1\in F\setminus U_1;$ suppose $x_1\in U_{n_1}.$ Then pick $x_2 \in
F\setminus (U_1\bigcup \cdots U_{n_1}).$ Suppose $x_k$ has been chosen and $%
x_k\in U_{n_k}.$ Choose $x_{k+1}\in F\setminus (U_1\bigcup \cdots U_{n_k}).$
These points $\{x_k\}$ must all be distinct. The sequence $\{x_k\}$ has a
subsequence that converges weak* sequentially to a point $y\in F.$ We assume
that $y\in \quad (some) \quad U_n.$

Now let $k^{\prime}$ be such that $n_{k^{\prime}}>n.$ Then $x_k\not\in U_n$
for all $k>k^{\prime}.$ Hence $U_n$ contains only finitely many points of $%
\{x_k\}$ and so the subsequence we found above cannot converge to $y$ in
weak* topology. This is a contradiction. \qed \newline

\begin{thm}[Banach-Alaoglu]
If $X$ is a normed space then $S*=\{x^{\prime}\in X^{\prime}:
||x^{\prime}||\leq 1\}$ is weak* compact.
\end{thm}

\begin{thm}
If $(X,T)$ is compact and if there exist continuous functions $%
\{f_{n}:X\rightarrow \mathbb{R}\}$ that separate points in $X$ (i.e. for any 
$x,y\in X$, $x\neq y,$ there is $n$ such that $f_{n}(x)\neq f(y)$), then $%
(X,T)$ is metrizable.
\end{thm}

\begin{thm}
If $X$ is a separable normed linear space and $K\subset X^{\prime}$ is weak*
compact, then $(K, W*)$ is metrizable.
\end{thm}

\textbf{Proof:}\quad By the above theorem we need only to find a countable
family of continuous functions from $(K,W\ast )$ to $\mathbb{R}$ which
separates points in $K.$ Let $x_{n}\in X$ and $\{x_{n}\}$ be dense in $X.$
Let $\Lambda _{n}:X^{\prime }\rightarrow C$ be defined as $\Lambda
_{n}(x^{\prime })=x^{\prime }(x_{n}).$ Then each $\Lambda _{n}$ is $W\ast $
continuous (by definition of weak* topology). Also $\{\Lambda _{n}\}$
separates points in $K.$ In fact, if $x^{\prime }\neq y^{\prime }$ are two
elements in $X^{\prime }$ and $\Lambda (x^{\prime })=\Lambda _{n}(y^{\prime
})$ for all $n,$ then $x^{\prime }$ and $y^{\prime }$ coincide on a dense
subset of $X$ and $x^{\prime }=y^{\prime }.$ A contradiction.

\begin{thm}[Weak* Compactness Theorem]
If $X$ is a separable normed linear space then the bounded sets in $%
X^{\prime}$ are weak* conditionally sequentially compact. That is, if $X$ is
separable and $x_n^{\prime}\in X^{\prime}$ with $||x^{\prime}_n||\leq A,$
then there is $x^{\prime}_0\in X^{\prime}$ with $||x^{\prime}_0||\leq A$ and
a subsequence $x^{\prime}_{n_k}$ such that $x^{\prime}_{n_k}\rightarrow
x^{\prime}_0$ in weak* topology, i.e., at each $x\in X,$ $%
x_{n_k}^{\prime}(x)\rightarrow x^{\prime}(x)$ as $k\rightarrow \infty.$ (cf.
page 22. Butzer)
\end{thm}

\textbf{Proof:}\quad The proof is obtained by putting together the
Banach-Alaoglu theorem and the above theorem.\qed\newline

\begin{cor}[Weak* Compactness Theorem for $L^p(\mathbb{R}),$ $1<p\leq
\infty. $]
For $1<p\leq \infty,$ if $||f_n||_p\leq A$ then there is $f_0,$ $%
||f_0||_p\leq A$ and $\{f_{n_k}\}$ so that for each $g\in L^{p^{\prime}},$ $%
\displaystyle{\int f_{n_k}g\rightarrow \int f_0g}.$
\end{cor}

\textbf{Proof:}\quad $L^p(\mathbb{R})$, $1<p\leq\infty,$ are conjugate
spaces of $L^{p^{\prime}}(\mathbb{R})$, $1\leq p<\infty,$ which are
separable. \qed

\begin{cor}
Let $\mu_n\in M$ (all finite Borel measures on $\mathbb{R}^n$) be such that $%
||\mu_n||_M\leq A$ for all $n.$ Then there is $n_k\rightarrow \infty$ and $%
\mu\in M$ so that $\mu_{n_k}\rightarrow \mu $ in weak* topology on M, that
is, for any $f\in C_0(\mathbb{R}^n),$ $\displaystyle{\int
fd\mu_{n_k}\rightarrow \int fd\mu}.$ Moreover, 
\begin{equation*}
||\mu||\leq \liminf_{k\rightarrow\infty}||u_{n_k}||.
\end{equation*}
\end{cor}

\textbf{Proof:}\quad $M(\mathbb{R}^n)$ is the conjugate space of $C_0(%
\mathbb{R}^n)$ ( by Riesz's theorem) which is separable. \qed \newline

\begin{cor}[Weak* Compactness Theorem for $L^1(\mathbb{R})$]
Let $f_n\in L^1(\mathbb{R})$ such that $||f_n||_1\leq K$ for all $n$. Then
there exist a subsequence $f_{n_k}$ and $\mu\in M(\mathbb{R})$ such that 
\begin{equation*}
\lim_{k\rightarrow\infty}\int_{\mathbb{R}}f_{n_k}(x)g(x)dx =\int_{\mathbb{R}%
} g(x)d\mu(x)
\end{equation*}
for each $g\in C_0(\mathbb{R}).$
\end{cor}

\textbf{Proof:}\quad We may view each $f_n$ as an element of $M(\mathbb{R}),$
if we identify $f_n(x)$ with $f_n(x)dx.$ Moreover, $||f_n(x)dx||_{M(\mathbb{R%
})}=||f_n||_1\leq K$ for all $n.$ The corollary follows. \qed \newline

\subsection{Dual or Conjugate Operators and Adjoint Operators}

Let $X,Y$ be Locally convex linear topological spaces. Let $T$ be a linear
operator on $D(T)\subset X$ into $Y.$ Let $\{x^{\prime},y^{\prime}\}$ be a
point in $X^{\prime}\times Y^{\prime}$ satisfying 
\begin{equation*}
<Tx,y^{\prime}>=<x, x^{\prime}>\,\,\forall x\in D(T).
\end{equation*}
Then $x^{\prime}$ is determined uniquely by $y^{\prime}$ iff $D(T)$ is dense
in $X.$

In this case, a linear operator $T^{\prime}$ defined by $T^{\prime}y^{%
\prime}=x^{\prime}$ is called the \emph{dual or conjugate operator} of $T.$
Its domain is the set of all $y^{\prime}\in Y^{\prime}$ such that there
exists $x^{\prime}\in X^{\prime}$ satisfying $<Tx,y^{\prime}>=<x,x^{\prime}>$
for all $x\in D(T).$

Let $X$ and $Y$ be complex Hilbert spaces. Let $E_{X}$ be the operator that
associates to each $y\in X$ the linear functional $x\rightarrow <x,y>.$ Then 
$E_{X}$ is a `conjugate-linear' isometry of $X$ onto $X.$ Let $T$ be a \emph{%
densely defined linear operator} from $X$ into $Y$. The \emph{adjoint}
operator of $T$ is the operator $T^{\ast }$ defined by 
\begin{equation*}
T^{\ast }=E_{X}^{-1}T^{\prime }E_{Y},
\end{equation*}%
where the domain of $T^{\ast }$ is the set of all $y$ for which $%
(E_{X}^{-1}T^{\prime }E_{Y})(y)$ is defined.

The notion of transposed conjugate matrix may be extended to the notion of
adjoint operator in Hilbert spaces. In contrast, the notion of transposed
matrix may be extended to the notion of dual operator in locally convex
linear topological spaces.

Clearly, 
\begin{equation*}
D(T^{\ast })=\{y:E_{Y}(y)\in D(T^{\prime })\}=\{y:x\rightarrow <Tx,y>\,\,%
\mbox{is continuous on $D(T)$}\}.
\end{equation*}
One can show that $y\in D(T^{\ast })$ if and only if there exists a $y^{\ast
}\in X$ such that 
\begin{equation*}
<Tx,y>=<x,y^{\ast }>
\end{equation*}%
holds for all $x\in D(T).$ In this case, $y^{\ast }=T^{\ast }y.$ If $T\in
L(X,Y),$ then $T^{\ast }\in L(Y,X)$ and $||T^{\ast }||=||T||.$ In general,
if $\overline{D(T)}=X,$ then $T^{\ast }$ is a closed linear operator.

It is known that $\overline{R(T)}=N(T^{\ast })^{\perp }.$ We include a
proof. If $y\in \overline{R(T)},$ then there exist $x_{n}\in D(T)$ such that 
$Tx_{n}\rightarrow y.$ Take $z\in N(T^{\ast }).$ Then $<Tx_{n},z>=<x_{n},T^{%
\ast }z>=0$ and so $<y,z>=0.$ This proves that $\overline{R(T)}\subset
N(T^{\ast })^{\perp }.$ To prove the opposite inclusion, we assume by
contradiction that there exists $p\in N(T^{\ast })^{\perp }$ but $p\not\in 
\overline{R(T)}.$ Then there is $f$ in the Hilbert space $X$ such that $%
<y,f>=0$ for all $y\in \overline{R(T)}$ and $<p,f>=1.$ Let $x\in D(T)$ and 
\emph{assume that $D(T)$ is dense in $X$}. Since $<x,T^{\ast }f>=<Tx,f>=0$
for all $x\in D(T),$ $T^{\ast }f=0.$ Since $p\in N(T^{\ast })^{\perp },$ $%
<p,y>=0$ for all $y$ with $T^{\ast }y=0.$ It follows that $<p,f>=0.$ This
contradiction proves that $\overline{R(T)}\supset N(T^{\ast })^{\perp }.$

Let $X$ and $Y$ be complex Hilbert spaces. An operator $U\in L(X,Y)$ is said
to be unitary if $U^{\ast }U=I_{X}$ (the identity on $X$) and $UU^{\ast
}=I_{Y}$ (the identity on $Y$). These two equations imply that $R(U)=Y,$ $%
D(U^{\ast })=Y,$ and $R(U^{\ast })=X.$ Given $U\in L(X,Y),$ the following
statements are equivalent: (1) $U$ is unitary; (2) $R(U)=Y$ and $U$
preserves the inner product; (3) $U$ is an isometric mapping of $X$ onto $Y.$

\section{Acknowlegments}

These notes were written by the first author in preparation for a series of
talks given on harmonic analysis through a succession of seminars at the
mathematics department of California State University Sacramento (CSUS).
Later on, the second author joined in a full collaborative effort to revise
and edit the notes and make them appropriate for publication as a graduate
level textbook. We are distinctly grateful to the faculty of the mathematics
department of CSUS for their helpful insights and support in preparation of
these notes. In writing the present manuscript, we would also like to
acknowledge that we were greatly inspired by Professor Henry Helson's
classic book, \emph{Harmonic Analysis}. Finally, we would like to express
our sincere appreciation to Professor Calixto Calderon of the University of
Illinois at Chicago for reviewing the final draft of the manuscript and
making helpful suggestions to improve it.

\section{Dedication}

The first author would like to sincerely express his gratitude to Mrs.
Zhenyan Zhou, his late wife, for her affectionate support and encouragement
during the writing of these notes. Without her unwavering and long time
support, the present work would not have been possible. The second author
also would like to acknowledge and express his gratefulness to Mrs. Mahin
Aliabadi Siadat, his late mother, for her never-ending encouragement and
loving support to persist in this collaboration, towards its successful
conclusion. We dedicate the present work to these highly honorable and
dedicated women. Although our loved ones are no longer with us, their
memories will for ever last in our hearts.

\end{document}